%% file: main_extended.tex
\documentclass[letterpaper, 10 pt, conference]{ieeeconf}  
\IEEEoverridecommandlockouts                              
\usepackage[utf8]{inputenc} 
\usepackage[T1]{fontenc}    
\usepackage{url}            
\usepackage{booktabs}       
\usepackage{amsmath,amssymb,amsfonts,mathrsfs}       
\usepackage{nicefrac}       
\usepackage{microtype}      
\usepackage{xcolor} 
\usepackage{xr-hyper,hyperref}

\usepackage{wasysym}
\usepackage{graphicx,float}
\usepackage{capt-of}
\usepackage{textcomp,comment}  
\let\labelindent\relax 
\usepackage{stfloats,empheq,epstopdf,enumitem}
\usepackage{tikz}
\usetikzlibrary{arrows,shapes,positioning,automata,backgrounds,petri}
\usepackage{graphics} 
\usepackage{pgfplots}
\usepackage{caption}
\newcommand{\cY}{\ensuremath{\mathscr{Y}}}
\newcommand{\cX}{\ensuremath{\mathscr{X}}}
\newcommand{\cZ}{\ensuremath{\mathscr{Z}}}
\newcommand{\cL}{\ensuremath{\mathscr{L}^2_{\rho_{\cX}}}}
\newcommand{\cH}{\ensuremath{\mathscr{H}}}
\newcommand{\cT}{\ensuremath{\mathscr{T}}}
\newcommand{\cK}{\ensuremath{\mathscr{K}}}
\newcommand{\Ss}{\ensuremath{\mathcal{S}_{\mathscr{X}}}}
\newcommand{\io}{\ensuremath{L_{\cK}}}
\newcommand{\lmin}{\ensuremath{\breve{\lambda}}}
\newcommand{\Var}{\ensuremath{\text{var}}}

\newtheorem{assumption}{Assumption}
\newtheorem{lemma}{Lemma}
\newtheorem{remark}{Remark}
\newtheorem{proposition}{Proposition}
\newtheorem{theorem}{Theorem}

\allowdisplaybreaks
\graphicspath{{Figures/}}

\title{\LARGE \bf
Error analysis of regularized trigonometric linear regression\\ with unbounded sampling: a statistical learning viewpoint
}

\author{Anna Scampicchio, Elena Arcari, Melanie N. Zeilinger
\thanks{The Authors are affiliated with the Institute of Dynamic Systems and Control, ETH Z\"urich.
        {\tt\small $\lbrace$ascampicc,earcari,mzeilinger$\rbrace$@ethz.ch}}%
}

\begin{document}

\maketitle
\thispagestyle{empty}
\pagestyle{empty}

\begin{abstract}
The effectiveness of non-parametric, kernel-based methods for function estimation comes at the price of high computational complexity, which hinders their applicability in adaptive, model-based control. Motivated by approximation techniques based on sparse spectrum Gaussian processes, 
we focus on models given by regularized trigonometric linear regression.  
This paper provides an analysis of the performance of such an estimation set-up within the statistical learning framework. 
In particular, we derive a novel bound for the sample error in finite-dimensional spaces, accounting for noise with potentially unbounded support. Next, we study the approximation error and discuss the bias-variance trade-off as a function of the regularization parameter by combining the two bounds. 
\end{abstract}

\section{INTRODUCTION}
Non-parametric approaches for regularized function estimation are a key tool in machine learning, and have been successfully applied to, e.g., system identification \cite{Pillonetto2014} and learning-based control~\cite{Hewing2020},~\cite{TK18}. Nevertheless, their applicability in real-time scenarios is hindered by their high computational complexity, which scales cubically with the number of data. 
The strategies proposed
to enable fast adaptation of kernel-based methods can be grouped into two main categories: input location selection, and low-rank approximations of the kernel \cite{scalableGP}. 
In this second class of approaches, a vast success was achieved by sparse spectrum Gaussian processes \cite{lazaro2010sparse,rudi2017}, 
where operations on the (shift-invariant) kernel yield a parametric approximation by means of linear combinations of Fourier features.\\ 
In this paper, we draw inspiration from the latter method and perform regularized regression within a finite-dimensional hypothesis space $\cH$ defined by a span of $E$ predefined trigonometric functions. 
Such a set-up relaxes the assumption of having shift-invariant kernels, and results more robust against potential basis function mis-specification compared to non-regularized, projection-based approaches \cite{hoerl1970}; for a review of parametric methods based on Fourier features, 
we refer to \cite[Chapter 1.7]{Tsybakov2009IntroductionTN}.  
Our goal is to assess the performance of the proposed estimator as a key step towards deriving 
reliable guarantees for data-driven, model-based control schemes that leverage such a model (see, e.g., \cite{pmlr-v70-pan17a,arcari2021,arcari22robot}).\\ 
We frame 
this analysis in the statistical learning set-up \cite{Cucker02,cuckerzhou2007}. 
The function to be estimated (i.e., the \textit{regression function} $f_{\rho}$) is defined as the minimizer of the expected risk over a (partially) unknown probability distribution, jointly defined over the input-output spaces, and from which i.i.d.~samples are drawn. Consequently, this formulation can also handle fully nonlinear measurement models. 
Furthermore, 
$f_{\rho}$ is generally assumed to belong to the space of square-integrable functions $\mathscr{L}^2$, and the hypothesis space is typically taken as an infinite-dimensional Reproducing Kernel Hilbert Space (RKHS), which is 
related to $\mathscr{L}^2$ by interpolation spaces arguments (\cite[Theorem 2]{Smale2007LearningTE}, \cite{lunardi2009}). Differently from classic non-parametric set-ups, the regression function is not assumed to belong to the hypothesis space. Thus, two objects can be therein defined: the actual data-based estimate $f_z$ and its data-free limit $f_{\cH}$. The goal of error analysis consists in quantifying the approximation error, or \textit{bias}, $\|f_{\rho} - f_{\cH}\|_{\cL}$, and the sample error, or \textit{variance}, $\|f_{\cH} - f_z\|_{\cL}$. As regards the latter, results abound in the statistical learning literature. Most of them deal with probability measures on the outputs that have bounded support, and thus obtain bounds leveraging concentration inequalities 
such as Hoeffding's or Bennett's \cite{Boucheron2004}, \cite[Chapter 3.1]{cuckerzhou2007}. Works in this direction are, e.g., \cite{Wu06learningrates,cuckerzhou2007,cucker2008bestchoices,mendelson2010,WANG201155}. Contributions considering unbounded sampling include \cite{Caponnetto2007OptimalRF,Wang2011,Guo2013ConcentrationEF}. The bounds therein derived leverage
the so-called moment hypothesis, which relaxes the boundedness assumption, and hold also for (sub-)Gaussian noises. Such results rely on the computation of covering numbers quantifying the capacity of the hypothesis space \cite{zhou2003}, and showcase optimal rates of convergence; nevertheless, they tend to be of limited practical relevance in the non-asymptotic case due to the large values of the multiplicative coefficients, which are often furthermore
difficult to compute.\\
In this work, we perform error analysis for finite-dimensional hypothesis spaces given by trigonometric functions. 
Our first contribution is a sample error bound, which 
is less conservative than the ones available in the literature even if it accounts also for noises with unbounded supports. 
Our second contribution consists in studying the bias-variance trade-off of the regularized trigonometric regression set-up. To this end, we obtain two bounds on the approximation error, combine them with the sample complexity result and analyze the conditions ensuring the existence of a unique value of the regularization parameter $\gamma$ returning the optimal trade-off. The differences of the two approaches for estimating the regularization parameter are investigated in a Monte Carlo study, which shows that one of the two criteria returns a value of $\gamma$ that captures the oracle behavior (i.e., minimizing the overall error),
thus leading to fast estimation schemes that do not need preliminary hyper-parameter selection. 

\section{PROBLEM SET-UP}\label{sec:problem}
Let the metric space of inputs $\cX$ be a compact subset of $\mathbb{R}$: without loss of generality, we take $\cX = [-X/2,\, X/2]$ for some $X \in \mathbb{R}_+$ (the scalar case is presented just for ease of visualization: the multi-dimensional is a straightforward extension). The output space is assumed to be $\cY = \mathbb{R}$. There is a probability measure $\rho$ defined over $\cZ = \cX \times \cY$ that decomposes into $\rho_{\cX}(x)$ and $\rho(y|x)$ according to Fubini's Theorem. In the considered 
set-up, the probability measure defined on $\cX$ is the standard uniform: denoting with $\mu$ the Lebesgue measure, we have that $\rho_{\cX}(A) = \mu(A \cap \cX)/\mu(\cX) = \mu(A \cap \cX)/X$ for any set $A$ in the $\sigma-$algebra of interest. In this way, $\rho_{\cX}$ is a Borel non-degenerate, $\sigma-$finite measure. As regards $\rho(y|\cdot)$, we assume it is unknown and defined over $\mathbb{R}$.\\
Having $N$ independent 
samples drawn from $\rho$ collected in the data-set $\mathcal{D} = \{(x_t,\,y_t)\}_{t=1}^N$, the goal is to estimate the regression function 
\begin{equation}
    f_{\rho}(x) = \int_{\cY}yd\rho(y|x). \label{eq:regressionfunction}
\end{equation}
We make use of the following Assumption.
\begin{assumption}
The regression function $f_{\rho}$ belongs to the space of square-integrable functions on $\cX$, denoted by $\cL$, and is such that $\|f_{\rho}\|_{\cL} = \sqrt{\int_{\cX} f^2(x)d\rho_{\cX}(x)} = \sqrt{\int_{\cX}f^2(x)d\mu(x)/X} \leq B_f$. Moreover, we also have that $\sigma^2_{\rho} = \int_{\cX} \sigma^2_{\rho}(x)d\rho_{\cX}(x) = \int_{\cZ} (y-f_{\rho}(x))^2d\rho \leq B_{\sigma}^2$. \hfill$\square$
\label{assume:bounds}
\end{assumption}
In other words, we assume to have access to bounds on the energy of the unknown function to be estimated, and on the variance of the additive noises.\\
The space $\cL$ is a separable Hilbert space whose complete orthonormal basis by means of trigonometric functions \cite{akhiezer2013theory} is given by 
\begin{align}
& \Bigg\lbrace \sqrt{2}\sin\Big(\frac{2\pi q}{X}x\Big), \sqrt{2}\cos\Big(\frac{2\pi q}{X}x\Big) \Bigg\rbrace_{q\in \mathbb{N}} \\&= \lbrace \bar{\varphi}_q^s(x),\, \bar{\varphi}_q^c(x) \rbrace_{q\in\mathbb{N}} \quad \text{with } x \in \cX.  
\label{eq:basis}
\end{align}
Accordingly, any function $f \in \cL$ can be expressed as $f(\cdot) = \sum_{q\in\mathbb{N}}\alpha^s_q\bar{\varphi}^s_q(\cdot) + \alpha^c_q\bar{\varphi}^c_q(\cdot)$, which will be also compactly written as $f(\cdot) = \sum_{q\in\mathbb{N}}\alpha_q \bar{\varphi}_q(\cdot)$, with $\sum_{q\in \mathbb{N}}\alpha_q^2 < \infty$. Within this representation, we denote the target function as $f_{\rho}(\cdot) = \sum_{q\in\mathbb{N}}\bar{\alpha}^s_q\bar{\varphi}^s_q(\cdot) + \bar{\alpha}^c_q\bar{\varphi}^c_q(\cdot) = \sum_{q\in\mathbb{N}}\bar{\alpha}_q \bar{\varphi}_q(\cdot)$.\\
Function estimation in $\cL$ cannot be performed, because pointwise evaluation is not well defined. Therefore, we perform such a task within a hypothesis space having the structure of a RKHS. Specifically, we consider the RKHS obtained from a subset of functions in \eqref{eq:basis} with cardinality~$E$, where $E$ is chosen according to our computational capacity. Denote by $Q$ the set of selected frequencies, i.e., $Q = \{q_j\}_{j=1}^{E/2} \subset \mathbb{N}$, and consider the following functions extracted from \eqref{eq:basis} using $Q$ defined, for $j=1,...,E/2$, as
\begin{equation}
\varphi_i(x) = 
    \begin{cases}
    \bar{\varphi}^s_{q_j}(x), \: i=j\\
    \bar{\varphi}^c_{q_j}(x), \: i = j+\frac{E}{2}.
    \end{cases}
    \label{eq:basisRKHS}
\end{equation}
Then, the RKHS of interest is the one induced by the following kernel:
\begin{align}
  \cK(x_a,x_b) = \phi^{\top}(x_a)\Sigma_{\alpha}\phi (x_b),\label{eq:ourkernel}
\end{align}
where $\Sigma_{\alpha} = \text{diag}(\lambda_1,...,\lambda_E)$ is a positive definite matrix, and the vector $\phi(\cdot) \in \mathbb{R}^E$ is such that $\phi^{\top}(x) =  [\varphi_1(x) \: \dots \: \varphi_E(x)].$ Clearly, \eqref{eq:ourkernel} is a Mercer kernel (\cite[(6), p.346]{aronszajn50reproducing}; it satisfies Mercer's condition $\int_{\cX}\int_{\cX}\cK(x,x^{\prime})^2d\rho_{\cX}(x)d\rho_{\cX}(x^{\prime}) = \sum_{i=1}^E \lambda_i^2$, and it is non-stationary if and only if $\lambda_i \neq \lambda_{i+E/2}$ for all $i=1,…,E/2$. Furthermore, using the argument in \cite[Chapter 4.3]{steinwartSVM}), it holds that 
\begin{align}
C_{\cK} &= \sup_{x_a,x_b \in \cX} \sqrt{\cK(x_a,x_b)} \notag \\ &\leq \sqrt{\sum_{i=1}^{E/2}\max\{\lambda_i,
, \lambda_{i+E/2}\}} < +\infty.
    \label{eq:boundCk}
\end{align}
Being a Mercer kernel, we have from Moore-Aronszajn Theorem \cite{aronszajn50reproducing} that $\cK$ is in one-to-one correspondence with the Hilbert space of functions $(\cH, \langle \cdot,\, \cdot \rangle_{\cH})$, which is 
\begin{equation}
\cH = \{f \in \cL: f(\cdot) = \phi^{\top}(\cdot)\alpha,\; \alpha \in \mathbb{R}^E\}
\label{eq:RKHS}
\end{equation}
with inner product given, for $f^{(\natural)}(\cdot) = \phi^{\top}(\cdot)\alpha^{(\natural)}$ and  $\natural= a ,\, b$: 
\begin{equation}
    \langle f^{(a)},\, f^{(b)} \rangle_{\cH} = \langle \Sigma_{\alpha}^{-1/2}\alpha^{(a)}, \Sigma_{\alpha}^{-1/2}\alpha^{(b)}\rangle_{2}. \label{eq:innerprodH}
\end{equation}
Within the hypothesis space, we can compute the estimate from the data-set $\mathcal{D}$ as follows. Consider the sampling operator $\Ss: \cH \rightarrow \mathbb{R}^N$ such that $\Ss(f) = [f(x_1) \ \dots \ f(x_N)]^{\top}$, together with its adjoint $\Ss^{\top}: \mathbb{R}^N \rightarrow \cH$ yielding $\Ss^{\top}c = \sum_{t=1}^N c_t\cK(x_t,\cdot)$. Thus, considering $Y = [y_1, ..., y_N]^{\top}$ and regularization parameter $\gamma>0$, we have
\begin{align}
    f_z &= \arg\min_{f \in \cH} \frac{1}{N}\sum_{t=1}^N (y_t - f(x_t))^2 + \gamma \|f\|_{\cH}^2 \label{eq:fz}\\
    &= \Big(\frac{1}{N}\Ss^{\top}\Ss + \gamma I \Big)^{-1}\frac{1}{N}\Ss^{\top}Y \label{eq:solSampOp}.
\end{align}
The aim of error analysis is to quantify the discrepancy between $f_z$ and $f_{\rho}$. To this end, we additionally consider the data-free limit of \eqref{eq:fz} as 
\begin{align}
    f_{\cH} &=  \arg\min_{f \in \cH} \int_{\cX}(f(x)-f_{\rho}(x))^2d\rho_{\cX}(x) + \gamma \|f\|_{\cH}^2 \label{eq:probdatafree}\\
    &= (\io + \gamma I)^{-1}\io f_{\rho}, \label{eq:soldatafree}
\end{align}
where $\io(f)(\bar{x})=\int_{\cX} \cK(\bar{x},x)f(x)d\rho_{\cX}(x)$ is an integral operator which, thanks to the properties of $\cK$, is (a) is self-adjoint and strictly positive, (b) continuous and compact, (c) satisties the Spectral Theorem \cite[Theorem 4.3]{cuckerzhou2007} with eigenpairs $\{(\varphi_i(\cdot), \lambda_i)\}_{i=1}^E$. Thanks to these properties, given an arbitrary $\cL$ function $f(x) = \sum_{q\in\mathbb{N}}\alpha_q\bar{\varphi}_q(x)$, using linearity and orthonormality of the basis, we have
\begin{equation}
L_{\cK}(f)(\bar{x}) = \sum_{i=1}^E \lambda_i \alpha_{i}^{\pi} \varphi_i(\bar{x}),
\label{eq:intopker}
\end{equation}
where we define the $i-$th component of the vector $\alpha^{\pi}$ for $i=1,...,E$, along the lines of \eqref{eq:basisRKHS}, as follows:
\begin{equation}
    \text{For }j=1,...,\frac{E}{2},\quad \alpha_i^{\pi} = 
    \begin{cases}
    \alpha_{q_j}^s, \qquad i=j\\
    \alpha_{q_j}^c, \qquad i=j+E/2.
    \end{cases}
    \label{eq:apiRKHS}
\end{equation}
Moreover, thanks to property (a), we can also define the $r$-th power of the integral operator\footnote{Note that the case $r=-1/2$ plays a crucial role in connecting the norms in $\cL$ and $\cH$ for functions in the hypothesis space. Indeed, considering $f(\cdot)= \sum_{i=1}^E \alpha_i\varphi_i(\cdot)$, one has by definition of $\cH$ that $\|f\|_{\cH}^2 = \|\Sigma_{\alpha}^{-1/2}\alpha\|^2_2 = \sum_{i=1}^E \alpha_i^2/\lambda_i$. On the other hand, we have that $\io^{-1/2}(f)(\cdot) = \sum_{i=1}^E \alpha_i/\sqrt{\lambda_i}\varphi_i(\cdot)$, and its $\cL$-norm is equal, by Parseval's Theorem, to $\sum_{i=1}^E \alpha_i^2/\lambda_i$. Therefore, we obtain that $\|f\|_{\cH}^2 = \|\io^{-1/2}f\|_{\cL}^2$.} as \cite{Cucker02}:
\begin{equation}
    L_{\cK}^r(f)(\bar{x}) = \sum_{i=1}^E \lambda_i^r \alpha_{i}^{\pi} \varphi_i(\bar{x}).
    \label{eq:powerIntOp}
\end{equation}

In the following, we study the sample error $\|f_z - f_{\cH}\|_{\cL}$ introduced by the finiteness of the data-set $\mathcal{D}$, and the approximation error $\|f_{\cH} - f_{\rho}\|_{\cL}$ determined by the choice of the hypothesis space. The two bound the overall error as $\|f_z - f_{\rho}\|_{\cL} \leq \|f_z - f_{\cH}\|_{\cL} + \|f_{\cH} - f_{\rho}\|_{\cL}$, which is to be minimized as a function of the regularization parameter~$\gamma$.

\section{SAMPLE ERROR}\label{sec:sampleerror}
In this Section we provide the novel result concerning the error between $f_z$ and $f_{\cH}$ introduced in \eqref{eq:solSampOp} and \eqref{eq:soldatafree}. Its proof can be found in Appendix \ref{sec:proofUniformBound}. \begin{theorem}
Let Assumption \ref{assume:bounds} hold. Consider $C_{\cK}$ as introduced in \eqref{eq:boundCk}, and define $\lmin = \min_{i=1,...,E} \lambda_i$. Then, with confidence at least $1-\delta$, it holds that
\begin{equation}
    \|f_z - f_{\cH}\|_{\cL} \leq \frac{C^3_{\cK}}{\gamma}\sqrt{\frac{B_f^2 + B_{\sigma}^2}{\lmin N\delta}}. \label{eq:uniformbound}
\end{equation}
\hfill$\square$
\label{prop:uniformbound}
\end{theorem}
\begin{remark}
We did not study bounds for $\mathbb{E}_{\cZ}[\rho_N(\|f_z - f_{\cH}\|_{\cL}]$, because they typically return conservative values. A result for unbounded sampling is given, e.g., in \cite[Proposition 20]{Lin2017distributed}. Note also that our probabilistic guarantees fall in the category of ``honest" bounds rather than ``exact" bounds, following the definitions given in \cite{Davies2007NonparametricRC}: this means that, for a user-chosen confidence level $\delta$, the result holds with confidence "at least $1-\delta$" and not with "exact probability $1-\delta$".
\end{remark}

\section{APPROXIMATION ERROR}\label{sec:approximationerror}
We now study the error due the choice of the hypothesis space $\cH$, i.e., the $\cL$-distance between the solution $f_{\cH}$ introduced in \eqref{eq:soldatafree} and the regression function $f_{\rho}$ defined in~\eqref{eq:regressionfunction}. We first provide an expression for $f_{\cH}$: letting the regression function be expressed through the basis functions of $\cL$ as $f_{\rho}(\cdot) = \sum_{q \in \mathbb{N}}\bar{\alpha}_q\bar{\varphi}_q(\cdot)$, and recalling the definition of the RKHS basis functions $\varphi_i(\cdot)$ in \eqref{eq:basisRKHS} and of the coefficients $\alpha_i^{\pi}$ in \eqref{eq:apiRKHS}, we have
\begin{equation}
    f_{\cH}(\cdot) = \sum_{i=1}^E \frac{\lambda_i}{\lambda_i + \gamma}\bar{\alpha}_i^{\pi}\varphi_i(\cdot). 
    \label{eq:lemmafH}
\end{equation}
Thanks to this result, we derive two bounds on the approximation error depending on different norms of the vector $\bar{\alpha}^{\pi}$ defined in \eqref{eq:apiRKHS}. The discussion of their performance is deferred to Section \ref{subsec:bv}. We present the result in the following Proposition, which is proven in Appendix \ref{sec:proofApproxErrs}. 
\begin{proposition}
In the trigonometric linear regression framework presented in Section \ref{sec:problem}, the approximation error $\|f_{\cH} - f_{\rho}\|_{\cL}$ admits the following upper bounds:
\begin{align}
    &\text{(a)}  \qquad \frac{\gamma}{\lmin + \gamma}\|\bar{\alpha}^{\pi}\|_2 + \sqrt{\sum_{q \in \mathbb{N}\setminus Q}\bar{\alpha}^2_q}  \label{eq:approxerror1}\\
    &\text{(b)}  \qquad  \|\bar{\alpha}^{\pi}\|_{\infty}\gamma\sum_{i=1}^E\frac{1}{\lambda_i} + \sqrt{\sum_{q \in \mathbb{N}\setminus Q}\bar{\alpha}^2_q}.\label{eq:approxerror2}
\end{align}
\hfill$\square$
\label{prop:approxerrs}
\end{proposition}

\section{BIAS-VARIANCE TRADE-OFF}\label{sec:biasvar}
In this section we combine the bounds on the sample and approximation errors derived in Theorem \ref{prop:uniformbound} and Proposition \ref{prop:approxerrs}, respectively, and study the estimated overall error~$\|f_z - f_{\rho}\|_{\cL}$ as a function of the regularization parameter~$\gamma$. We perform our analysis after the RKHS $\cH$ has been completely specified, i.e., after having fixed $Q$ and~$\{\lambda_i\}_{i=1}^E$. \\
The main result is presented in the following Proposition proven in Appendix \ref{sec:proofTradeOff}.
\begin{proposition} 
\begin{enumerate}[label=(\alph*)]
\item[]
\item Consider the approximation error bound as in \eqref{eq:approxerror1}. Then, if the number of data $N$ and the confidence parameter $\delta$ are such that
\begin{equation}
    \sqrt{N\delta} > \frac{C_{\cK}^3}{\lmin^{3/2}}\sqrt{\frac{B_f^2 + B_{\sigma}^2}{\sum_{i=1}^E (\bar\alpha_i^{\pi})^2}},
    \label{eq:bvcond}
\end{equation}
there exists a unique $\gamma=\hat{\gamma}_{(a)}$ minimizing the estimated error $\|f_z - f_{\rho}\|_{\cL}$.\\
\item Take now the approximation error bound as in \eqref{eq:approxerror2}. Then, there always exist a unique $\gamma = \hat{\gamma}_{(b)}$ minimizing the estimated error $\|f_z - f_{\rho}\|_{\cL}$. \hfill$\square$
\end{enumerate}

\label{prop:tradeoff}
\end{proposition}
The closed-form expressions for $\hat{\gamma}_{(a)}$ and $\hat{\gamma}_{(b)}$ are provided in the proof.

\section{DISCUSSION}\label{sec:discussion}
We first study the performance of the sample error bound provided in Section \ref{sec:sampleerror} by comparing it with other bounds given in \cite{Smale2007LearningTE} and \cite{Lin2017distributed}. Next, we discuss the result of Proposition \ref{prop:tradeoff}, especially showcasing the capability of $\gamma^{(b)}$ to capture the behaviour of the oracle $\gamma$ minimizing the overall error. 

\subsection{Comparison with sample error bound in \cite[Theorem 5]{Smale2007LearningTE}}\label{sec:discSZ}
In the numerical set-up we assume that a uniformly distributed random noise with a Signal-to-Noise Ratio (SNR) of 150 affects the measurements of the regression function $f_{\rho}(x) = \sum_{q\in\mathbb{N}}\bar{\varphi}_q(x)\bar{\alpha}_q$ with $x \in [-1250,1250]$. Such a function is assumed to be characterized by 20 sine/cosine couples $\{\bar{\varphi}_q\}$, where $q$ is randomly drawn without repetitions from the set $\{1,...,30\}$. The hypothesis space $\cH$ is characterized by a subset of $E/2 = 10$ sine/cosine couples randomly selected among those that define the regression function.\\
We perform a Monte Carlo study of 500 runs, where at each step we draw a new set of basis functions defining the regression function and the hypothesis space.  Coefficients~$\bar\alpha_q$ of the regression function are drawn from a Gaussian distribution $\mathcal{N}(0,\lambda)$, where $\lambda$ is sampled from a uniform distribution on $(0,5)$, and also enters the definition of the hypothesis space as in \eqref{eq:innerprodH} as $\lambda_i = \lambda$ for all $i=1,...,E$. At each run, the number of data-points $N$ is randomly sampled from the set $\{100,101,...,1000\}$. We consider a confidence level of $\delta=0.1$. Then, we evaluate the sample error bounds corresponding to the minimum value of $\gamma$ satisfying the bound in \cite[Theorem 5]{Smale2007LearningTE}, and evaluate their relative difference with respect to the true sample error attained with such a $\gamma$. The results are displayed in Figure \ref{fig:SZandLGZboxplotBounds}. 
Both bounds decay as $1/\sqrt{N}$, but~\eqref{eq:uniformbound} evidences a more favorable behaviour in terms of the confidence level, at least for values of $\delta$ smaller than the solution of $1/\sqrt{\delta} = \log(4/\delta)$ in $(0,1]$, that is $\approx 0.0539$. Conservatism in the bound in \cite[Theorem 5]{Smale2007LearningTE} is mostly due to the linear dependence on the output values bound, $M$.  
The explicit condition on $M$ ensuring bound \eqref{eq:uniformbound} to be more conservative is the following:
\begin{equation}
    M \leq \frac{C_{\cK}^2}{12}\sqrt{\frac{B_f^2 + B_{\sigma}^2}{\lmin \gamma}}\frac{1}{\sqrt{\delta}\log(4/\delta)}. \label{eq:boundM}
\end{equation}
Such a value tends to be very small: e.g., in the Monte Carlo test, the bound \eqref{eq:boundM} returned a mean value of $3.50 \pm 2.02$, while the true value $M$ emerging from the (quite favorable) SNR attained a mean value of $39.02 \pm 14.66$.

\subsection{Comparison with sample error bound in \cite[Proposition 20]{Lin2017distributed}}\label{sec:discLGZ}
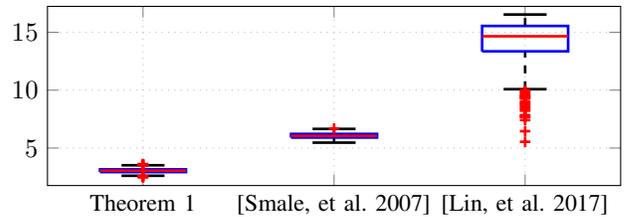
\begin{figure}[h]
    \centering
    \input{SZandLGZcomparison.tex}
    \caption{Behaviour of the sample error bounds in the Monte Carlo trials in Sections \ref{sec:discSZ} and \ref{sec:discLGZ}. The adopted score is the difference between bound and true sample error, normalized by the true sample error. For Theorem \ref{prop:uniformbound}, such an error attains a mean value of $21.44 \pm 4.093$, while for the bound in \cite{Smale2007LearningTE} it is $440.03 \pm 99.33$, and $3.40\times 10^6\pm 3.32 \times 10^6$ for the one in \cite{Lin2017distributed}. We display the values in logarithmic scale to facilitate visualization.  
   }
    \label{fig:SZandLGZboxplotBounds}
\end{figure}

We consider the same numerical set-up as the previous section, and we translate the bound of \cite[Proposition 20]{Lin2017distributed} into a statement of the same type as Theorem \ref{prop:uniformbound} by using Markov's inequality. To further adapt to the context given in Section \ref{sec:problem}, we set $p=2$ and~$\mathcal{N}(\gamma) = \sum_{i=1}^E \lambda_i/(\lambda_i + \gamma)$. The bound of \cite[Proposition 20]{Lin2017distributed} shows a slower behaviour in the number of data $N$ with respect to \eqref{eq:uniformbound}; moreover, it depends on the approximation error, which is generally not known. We performed the Monte Carlo study by setting~$\|f_{\cH} - f_{\rho}\|_{\cL}$ to its true value, and the results are very conservative, as displayed in Figure \ref{fig:SZandLGZboxplotBounds}.

\subsection{On the choice of $\gamma$ in view of the bias-variance trade-off} \label{subsec:bv}
We now perform a Monte Carlo study to discuss the results given in Proposition \ref{prop:tradeoff}. Consider $\cX = [-5\times 10^5,\,5\times 10^5]$ as input domain. The regression function is characterized by 30 sine/cosine pairs $\{\bar{\varphi}\}_{q=1}^{30}$, where each $q$ is randomly selected without repetitions from the set $\{1,...,100\}$, and each component of the vector of coefficients $\bar{\alpha}$ is drawn from a Gaussian random variable with zero mean and variance $\lambda=1$. The latter hyper-parameter also enters the definition of the RKHS $\cH$. The set of frequencies $Q$ is selected as a random subset with cardinality 10 from the set of those characterising the regression function. Fixing an SNR equal to 50, we draw 50 random regression functions and select the basis functions for the hypothesis space. The number of input/output pairs for each run is $N=2500$, and we consider a confidence parameter $\delta=0.5$. For each run, we compute $\gamma^{(a)}$ and $\gamma^{(b)}$ as in Proposition~\ref{prop:tradeoff}, compute the sample- and approximation error bounds as in Theorem~\ref{prop:uniformbound} and Proposition~\ref{prop:approxerrs}, and compare their values to the true errors yielded by $\gamma^{(a)}$ and $\gamma^{(b)}$. 
We observe that the bounds computed with $\gamma^{(a)}$ are closer to the true values. We display the values in Table \ref{tab:bv}.

\begin{table}[!h]
    \centering
    \begin{tabular}{|c|c|c|}
    \hline
        $\gamma_{(a)}$ & True value & Bound\\
         \hline
         $\|f_z - f_{\cH}\|_{\cL}$ &$0.036 \pm 0.007$ &$0.419 \pm 0.033$  \\
         \hline
         $\|f_{\cH} - f_{\rho}\|_{\cL}$ \tiny{(a)} &$9.313 \pm 0.077$ &$10.05 \pm 0.992$ \\
         \hline
    \hline
         $\gamma_{(b)}$ & True value & Bound\\
         \hline
         $\|f_z - f_{\cH}\|_{\cL}$  & $0.387 \pm 0.076$ & $14.04 \pm 1.942$ \\
         \hline
         $\|f_{\cH} - f_{\rho}\|_{\cL}$ \tiny{(b)} &$7.351 \pm 0.575$ & $20.41 \pm 2.193$\\
         \hline
    \end{tabular}
    \vspace{0.5em}
    \caption{Overall values (mean $\pm$ standard deviation) of sample and approximation error bounds compared to the true errors.}
    \label{tab:bv}
\end{table} 
The test above described was performed fixing the regularization parameter, and focused on the single errors. If we instead consider the overall error $\|f_z - f_{\rho}\|_{\cL}$, and compare values of $\gamma^{(a)}$ and $\gamma^{(b)}$ with the oracle value $\gamma^*$ (obtained via grid search) minimizing it, we observe that $\gamma^{(b)}$ is the one that performs best. The poor performance of $\gamma^{(a)}$ is due to the fact that the condition in~\eqref{eq:bvcond} needs a large number of data to be satisfied, and this leads to an overestimation of the regularization parameter. In this specific test, $\gamma^*$ was located at the minimum value of the grid, i.e. $\gamma^*=0.1$; the mean values for $\gamma^{(a)}$ and $\gamma^{(b)}$ were $7.7703 \pm 0.2115$ and $0.2308 \pm 0.0230$, respectively. 

\section{CONCLUSIONS}
In this paper, we analysed the estimation errors occurring in regularized trigonometric regression within the statistical learning set-up. To the best of the Authors' knowledge, such a study was missing in the literature, that mostly focused on non-parametric methods or non-regularized trigonometric regression. We derived a novel bound on the sample error that does not require the support of the output distribution to be finite; numerical tests showed it to be less conservative than classical bounds, at least in the non-asymptotic regime. Next, we computed two bounds for the approximation error, and combined them with the sample error bound to retrieve a selection criterion for the regularization parameter~$\gamma$, optimizing the trade-off between estimated bias and variance. In particular, we showed that one of the two criteria yields a value of the regularization parameter that is close to the oracle, and thus can in principle be used to speed up hyper-parameter selection. We stress that such an analysis can be extended to any other orthogonal basis of $\cL$. Moreover, we foresee that the generality of such a set-up can have an impact on an abstract treatment of bias learning, which is a planned extension of the present work. Forthcoming research will also focus on asymptotic behavior in terms of number of data $N$ and of the basis functions $E$.

\bibliographystyle{IEEEtran}
\bibliography{Main}
 \newpage
\onecolumn
\appendix
We provide in Appendix 
\ref{sec:introstatlear} an introduction to error analysis in the statistical learning framework. Appendices 
\ref{sec:proofUniformBound}, \ref{sec:proofApproxErrs} and \ref{sec:proofTradeOff} report the proofs for the theoretical results stated in Sections \ref{sec:sampleerror}, \ref{sec:approximationerror} and \ref{sec:biasvar}, respectively. For ease of referencing, the statements of the Theorems used as benchmarks in Section \ref{sec:discussion} are given in Appendix 
\ref{sec:statementsdisc}. Appendix 
\ref{sec:discGZ} presents an additional comparison with the sample error bound with unbounded noise support reported in \cite{Wang2011}. In Appendix 
\ref{sec:addendaBV}, we show additional plots related to the experiment of Section \ref{subsec:bv}. A discussion on the benefits of regularization is given in Appendix 
\ref{subsec:regularization}.

\subsection{Statistical learning framework}\label{sec:introstatlear}

Let the metric space of inputs $\cX$ be compact, and let $\cY=\mathbb{R}$ be the space of outputs. There is a probability measure $\rho$ defined over $\cZ = \cX \times \cY$ that decomposes as $\rho(y|x)$ and $\rho_{\cX}(x)$ according to Fubini's Theorem: 
given an integrable function $\psi: \cZ \rightarrow \mathbb{R}$, it holds  
\begin{equation}
   \int_{\cZ} \psi(z)d\rho(z) = \int_{\cX}\Big(\int_{\cY}\psi(x,y)d\rho(y|x)\Big)d\rho_{\cX}(x). \notag
    \label{eq:fubini}
\end{equation}   
Assume to collect $N$ independent samples drawn from $\rho$ in the data-set $\mathcal{D} = \{(x_t,\,y_t)\}_{t=1}^N$. The aim of statistical learning is that of estimating the \textit{regression function} of $\rho$ defined as
\begin{equation}
    f_{\rho}(x) = \int_{\cY}yd\rho(y|x). \notag \label{eq:regressionfunction1}
\end{equation}
Such a function is the minimizer of the expected risk 
$\mathcal{I}[f] = \mathbb{E}_{\rho}[(y-f(x))^2] = \int_{\cZ} (y-f_{\rho}(x) + f_{\rho}(x) - f(x))^2d\rho = \mathcal{I}[f_{\rho}] + \int_{\cZ} ( f_{\rho}(x) - f(x))^2d\rho$ \cite{Niyogi99} and can be viewed as the first statistical moment of $\rho(y|x)$. Its variance is defined as
\begin{equation}
\sigma^2_{\rho}(x) = \int_{\cY}(y-f_{\rho}(x))^2d\rho(y|x), \notag
\label{eq:varRhox1}
\end{equation}
whose integral over the $\cX$ domain is 
\begin{equation}
\sigma^2_{\rho} = \int_{\cX} \sigma^2_{\rho}(x)d\rho_{\cX}(x) = \int_{\cZ} (y-f_{\rho}(x))^2d\rho = \mathcal{I}[f_{\rho}]. \notag
\label{eq:varRho1}
\end{equation}
Note that $\sigma^2_{\rho}$ represents the unavoidable cost in the minimization of $\mathcal{I}(f)$ and it measures how well conditioned $\rho$ is: in other words, it is analogous to the notion of condition number in numerical linear algebra \cite{Cucker02}.\\

The regression function is assumed to belong to a certain \textit{target space} $\cT$, but is not computable in practice because $\rho$ is not known. Therefore, the estimate $f_z$ is searched within a \textit{hypothesis space} $\cH$ that is amenable to perform computations. To this aim, $\cH$ is chosen as a Reproducing Kernel Hilbert Space (RKHS) with inner product $\langle \cdot, \, \cdot \rangle_{\cH}$, and the regression function is estimated by solving the following Tikhonov regularization problem:
\begin{equation}
    f_z = \arg\min_{f \in \cH} \frac{1}{N}\sum_{t=1}^N (y_t - f(x_t))^2 + \gamma \|f\|_{\cH}^2. \notag
    \label{eq:fz1}
\end{equation}
The key feature in RKHSs is that function evaluation at any point in the domain $\cX$ is well defined by a functional that is linear and continuous. From Moore-Aronszajn Theorem \cite{aronszajn50reproducing}, it results that the RKHS is in one-to-one correspondence with a positive semi-definite kernel operator (see, e.g., \cite[Definition 2.8]{cuckerzhou2007})
\begin{equation}
\cK:\cX \times \cX \rightarrow \mathbb{R} \notag
\label{eq:kernel1}
\end{equation}
such that the \textit{reproducing property} holds, i.e., $f(x)=\langle \cK(x,\cdot), f\rangle_{\cH}.$ These facts, together with Riesz-Frechet theorem (see, e.g., \cite{berberian1961introduction}[Chapter V, Theorem 1] 
) lead to the so-called Representer Theorem \cite{wahba2019representer}, stating that the solution $f_z$ is a linear combination of $\{\mathscr{K}(x_{t},\cdot)\}_{t=1}^{N}$, i.e., of the kernel sections centered at the given input locations in $\mathscr{X}$. 
The same result can be expressed via the \textit{sampling operator} $\Ss$ \cite{smaleshannonII,Smale2007LearningTE}:
\begin{lemma}
Let $\Ss: \cH \rightarrow \mathbb{R}^N$ be an operator such that $\Ss(f) = [f(x_1) \ \dots \ f(x_N)]^{\top}$, and consider its adjoint $\Ss^{\top}: \mathbb{R}^N \rightarrow \cH$ yielding $\Ss^{\top}c = \sum_{t=1}^N c_t\cK(x_t,\cdot)$. Introducing $Y = [y_1 \ \dots \ \textbf{}y_N]^{\top}$, we have
\begin{equation}
f_{z} =  \Big(\frac{1}{N}\Ss^{\top}\Ss + \gamma I \Big)^{-1}\frac{1}{N}\Ss^{\top}Y. \notag
\label{eq:solSampOp1}
\end{equation}
\label{lemma:fz}
\end{lemma}
\textit{Proof.}
Let us begin with the derivation of the expression for $\Ss^{\top}$. By definition of adjoint operator, one must have that $\langle \Ss f, c \rangle_{2} = \langle f, \Ss^{\top}c \rangle_{\cH}$.  Now,  $\langle \Ss f, c \rangle_{2} = \sum_{t=1}^N c_tf(x_t)$, which by the reproducing property is equal to $\sum_{t=1}^N c_t \langle f, \cK(x_t,\cdot)\rangle_{\cH} = \langle f, \sum_{t=1}^N c_t\cK(x_t,\cdot)\rangle_{\cH}$.  By inspection of the definition of adjoint operator, it follows that $\Ss^{\top}c =\sum_{t=1}^N c_t\cK(x_t,\cdot)$. \\

Let us now retrieve the expression for $f_z$. 
Write the objective as 
\begin{align*}
    &\frac{1}{N}\Big(\langle \Ss f, \Ss f \rangle_{2} + \|Y\|_{2}^2 - 2\langle y, \Ss f \rangle_{2}\Big) + \gamma \langle f,f\rangle_{\cH}\\
    =&\Big\langle \Big(\frac{1}{N}\Ss^{\top}\Ss + \gamma I\Big)f,f\Big\rangle_{\cH} - \frac{2}{N}\langle \Ss^{\top}y, f\rangle_{\cH} + \frac{1}{N}\|Y\|_{2}^2.
\end{align*}
Solution follows by taking the functional derivative of the last expression. \hfill$\blacksquare$

The primary question of interest in error analysis is about quantifying how well $f_z$ approximates $f_{\rho}$. To this aim, we introduce the data-free limit of $f_z$ as 
\begin{equation}
    f_{\cH} = \arg\min_{f \in \cH} \int_{\cX}(f(x)-f_{\rho}(x))^2d\rho_{\cX}(x) + \gamma \|f\|_{\cH}^2.\notag
    \label{eq:probdatafree1}
\end{equation}
Its solution is given by means of the \textit{integral operator} 
\begin{equation}\label{eq:intOp}
\io(f)(\bar{x}) = \int_{\cX} \cK(\bar{x},x)f(x)d\rho_{\cX}(x), \notag
\end{equation}
and its expression is (\cite{Cucker02}, II.2, Theorem 3; and III.6, Proposition 7) \begin{equation}
f_{\cH} = (\io + \gamma I)^{-1}\io f_{\rho}. \notag
\label{eq:soldatafree1}
\end{equation}
Having defined $f_{\rho}$, $f_z$ and $f_{\cH}$, consider the metric $\|f_z - f_{\rho}\|_{\natural}$, where $\natural$ indicates the type of norm of interest. The overall error thus decomposes as
\begin{equation}
    \|f_z - f_{\rho}\|_{\natural} \leq \|f_z - f_{\cH}\|_{\natural} + \|f_{\cH} - f_{\rho}\|_{\natural}. \notag
    \label{eq:overallerror}
\end{equation}
The first addendum is named \textit{sample error}, and indicates the error within the RKHS $\cH$ due to the fact that we are operating with a finite amount of data. The second is called \textit{approximation error} and arises from the choice of the hypothesis space. When considering the error in the space of square-integrable functions, the two terms are also called \textit{variance} and \textit{bias}, respectively. The choice on the size of $\cH$ has an opposite effect on them: the larger the hypothesis space is, the smaller the distance from $f_{\cH}$ to $f_{\rho}$ can be; on the other hand, the more complex the model is, the more data are required to fit it. 

\begin{remark}
The statistical learning viewpoint was originally juxtaposed to the so-called \textit{sampling theory} approach for function estimation \cite{smaleshannonI,Niyogi99}. The first can be viewed as a more flexible framework to perform error analysis, and comprises the latter as a special case. Indeed, one could envisage $f_{\rho}$ as the ``true" function to be estimated, assuming data are generated as $y_t = f_{\rho}(x_t) + e_t$ and having noises with zero mean and variance $\sigma_{\rho}^2(x_t)$. For further comments, please refer to \cite[Section 7]{smaleshannonII}. 
\end{remark}

\subsection{Proof of Theorem \ref{prop:uniformbound}}\label{sec:proofUniformBound}
Defining $\xi_t : \cZ \rightarrow \cH$ such that $\xi_t(\cdot) = (y_t - f_{\cH}(x_t))\cK(x_t,\,\cdot)$, it holds that $\mathbb{E}_{\cZ}[\xi_t](\cdot) = L_{\cK}(f_{\rho} - f_{\cH})(\cdot) = \gamma f_{\cH}(\cdot)$. From this, and recalling the definition of sampling operator, it follows that $f_z(x) - f_{\cH}(x)$ is equal to \cite{Smale2007LearningTE}
\begin{equation}
   \Big(\frac{1}{N}\Ss^{\top}\Ss + \gamma I \Big)^{-1}\Bigg[\frac{1}{N}\sum_{t=1}^N \xi_t(x) - \mathbb{E}_{\cZ}[\xi](x)\Bigg].\notag
\end{equation}
We can now study the $\cL-$ norm of the expression above. Since $\cX$ is compact and the measure $\rho_{\cX}$ on it defined is a probability measure, $\|f\|_{\cL} \leq \|f\|_{\infty}$ for any function $f \in \cL$: therefore, $\|f_z - f_{\cH} \|_{\cL}$ is upper bounded by
\begin{equation}
  \Big\|\Big(\frac{1}{N}\Ss^{\top}\Ss + \gamma I \Big)^{-1}\Big\|_{\infty}\Big\|\frac{1}{N}\sum_{t=1}^n\xi_t - \mathbb{E}_{\cZ}[\xi]\Big\|_{\infty}.  \notag
\end{equation}
Since the operator norm can be bounded by $\frac{C_{\cK}}{\gamma\sqrt{\lmin}}$ (the proof is reported at the end of this subsection), we can now study an upper bound for $\rho_{N}(\|f_z - f_{\cH}\|_{\cL} > \epsilon)$ which, for an arbitrary $\epsilon >0$, is
\begin{equation}
    \rho_{N}\Bigg(\Big\|\frac{1}{N}\sum_{t=1}^n\xi_t - \mathbb{E}_{\cZ}[\xi]\Big\|_{\infty}   > \frac{\epsilon\gamma \sqrt{\lmin}}{C_{\cK}} \Bigg).
    \label{eq:boundsrhoN}
\end{equation}
At an arbitrary input location $x \in \cX$ and a given $\bar{\epsilon} \in (0,1)$, Chebychev's inequality yields
\begin{equation}
    \rho_N\Big(\Big|\frac{1}{N}\sum_{t=1}^N \xi_t(x) - \mathbb{E}_{\cZ}[\xi](x)\Big|  > \bar{\epsilon}\Big) \leq
    \frac{\Var(\xi)(x)}{N\bar{\epsilon}^2}, \notag
\end{equation}
noting that $\{\xi_t\}_{t=1}^N$ are independent and identically distributed. Using this result, we can further bound \eqref{eq:boundsrhoN} as
\begin{equation}
    \rho_{N}(\|f_z - f_{\cH}\|_{\cL} > \epsilon) \leq \frac{C_{\cK}^2}{\gamma^2 \lmin}\frac{\|\Var(\xi)\|_{\infty}}{N\epsilon^2}.
    \label{eq:chebychevour}
\end{equation}
The variance term can be bounded as 
\begin{align*}
    \sup_{\bar{x}\in \cX} \Var(\xi)(\bar{x}) &\leq  \sup_{\bar{x}\in \cX} \int_{\cZ}\cK(\bar{x},x)^2(y - f_{\cH}(x))^2d\rho \leq C_{\cK}^4\int_{\cZ}(y - f_{\cH}(x))^2d\rho \leq B_f^2 + B_{\sigma}^2, \notag
\end{align*}
where the last inequality follows from the fact that $\int_{\cZ}(f(x)-y)^2d\rho - \int_{\cZ}(f_{\rho}(x)-y)^2d\rho = \|f - f_{\rho}\|^2_{\cL}$ for any $f:\cX \rightarrow \cY$ \cite{Smale2007LearningTE}, and that $\|f_{\cH} - f_{\rho}\|^2_{\cL} + \gamma\|f_{\cH}\|^2_{\cH} = \mathcal{J}(f_{\cH}) \leq \mathcal{J}(0)   = \|f_{\rho}\|^2_{\cL} \leq B_f^2$.
Coming back to \eqref{eq:chebychevour}, we have that 
\begin{equation}
    \rho_{N}\Big(\|f_z - f_{\cH}\|_{\cL} > \epsilon\Big) \leq \frac{C_{\cK}^6}{\gamma^2 \lmin}\frac{(B_f^2 + B_{\sigma}^2)}{N\epsilon^2} = \delta.
    \label{eq:lastep}
\end{equation}
The proof is concluded by retrieving the expression for $\epsilon$ from $\delta$ in the equality \eqref{eq:lastep}. \hfill$\blacksquare$
\paragraph*{Proof for operator norm bound} By definition, we look for a constant $\mathfrak{C}_{\infty}$ is such that, for any $u \in \cH$,  $\|(\Ss^{\top}\Ss/N + \gamma I)^{-1}u\|_{\infty} \leq \mathfrak{C}_{\infty}\|u\|_{\infty}$. By the reproducing property, $\Big\|\Big(\frac{1}{N}\Ss^{\top}\Ss + \gamma I \Big)^{-1} u\Big\|_{\infty} = \sup_{\bar{x} \in \cX} \Big| \Big\langle \Big(\frac{1}{N}\Ss^{\top}\Ss + \gamma I \Big)^{-1} u(\cdot), \cK(\bar{x},\cdot)  \Big\rangle_{\cH} \Big|$, which is further upper bounded by $C_{\cK}\Big\|\Big(\frac{1}{N}\Ss^{\top}\Ss + \gamma I \Big)^{-1}\Big\|_{\cH}\|u\|_{\cH}$ by Cauchy-Schwartz inequality and \eqref{eq:boundCk}. Now, by the bound on the operator norm in $\cH$ provided in \cite[Equation 3.5]{Smale2007LearningTE}, we have $\Big\|\Big(\frac{1}{N}\Ss^{\top}\Ss + \gamma I \Big)^{-1} u\Big\|_{\infty} \leq \frac{C_{\cK}}{\gamma}\|\io^{-1/2}\|_{\cL}\|u\|_{\cL} \leq \frac{C_{\cK}}{\gamma}\|\io^{-1/2}\|_{\cL}\|u\|_{\infty}$. The proof is concluded by deriving the operator norm for $\|\io^{-1/2}\|_{\cL}$, which is $\|\io^{-1/2}\|_{\cL} \leq 1/\sqrt{\lmin}$ because, for an arbitrary $f \in \cL$, $\|\io^{-1/2} f\|_{\cL} = \sqrt{\sum_{i=1}^E \frac{\alpha_i^2}{\lambda_i}} \leq \sqrt{\frac{1}{\lmin}\sum_{i=1}^E \alpha_i^2} \leq \frac{1}{\sqrt{\lmin}}\|f\|_{\cL}$. 

\subsection{Proof of Proposition \ref{prop:approxerrs}}\label{sec:proofApproxErrs}
Expressing the regression function as $f_{\rho} = \sum_{q\in\mathbb{N}}\bar{\alpha}_q\bar{\varphi}_q$ and $f_{\cH}$ as in \eqref{eq:lemmafH}, we apply the triangle inequality and Parseval's Theorem on $\|f_{\cH} - f_{\rho}\|_{\cL}$ and obtain
\begin{align}
    \|f_{\cH} - f_{\rho}\|_{\cL} &= \Big\|\sum_{i=1}^E \frac{\lambda_i}{\lambda_i + \gamma}\bar{\alpha}^{\pi}_{i}\varphi_i - \sum_{q \in \mathbb{N}}\bar{\alpha}_q\bar{\varphi}_q \Big\|_{\cL} 
    \leq \sqrt{\sum_{i=1}^E \Bigg(\frac{\gamma}{\lambda_i + \gamma}\Bigg)^2(\bar{\alpha}_{i}^{\pi})^2} + \sqrt{\sum_{q \in \mathbb{N}\setminus Q}\bar{\alpha}_{q}^2}. \notag
\end{align}
Let us now focus on the first term on the right-hand side. The first bound \eqref{eq:approxerror1} is obtained by considering $(\lambda_i + \gamma)^{-1} \leq (\lmin + \gamma)^{-1}$. As for the second, we take $\bar{\alpha}^{\pi}_{i} \leq \|\bar{\alpha}^{\pi}\|_{\infty}$, bound the square root of the sum as the sum of the square roots, and take $(\lambda_i + \gamma)^{-1} \leq (\lambda_i)^{-1}$.

\subsection{Proof of Proposition \ref{prop:tradeoff}}\label{sec:proofTradeOff}
(a) Consider the sample and approximation errors as obtained in \eqref{eq:uniformbound} and \eqref{eq:approxerror1}, respectively. 
Introducing the following notation:
\begin{equation}
A = C_{\cK}^3\sqrt{\frac{B_f^2 + B_{\sigma}^2}{N\delta\lmin}}, \quad b = \lmin, \quad
        B = \sqrt{\sum_{i=1}^E (\bar{\alpha}^{\pi}_{i})^2}, \quad C = \sqrt{\sum_{q \in \mathbb{N}\setminus Q}\bar\alpha_q^2},
\label{eq:Abbv}
\end{equation}
we have that the overall error can be bounded as follows:
\begin{equation}
    \|f_z - f_{\rho}\|_{\cL} \leq \frac{A}{\gamma} + \frac{B\gamma}{b + \gamma} + C = F(\gamma).
    \label{eq:symbBound}
\end{equation}
The function $F(\gamma)$ is always positive for $\gamma > 0$.
We aim at finding the condition for which there exists a unique, finite value of $\gamma$ minimizing $F(\gamma)$. To this end, let us study the first derivative:
\begin{equation}
   \frac{dF}{d\gamma} = 0 \longrightarrow \gamma^2(Bb-A) - 2Ab\gamma - Ab^2 = 0. \label{eq:firstDerBV}
\end{equation}
By applying Descartes' rule, we obtain that the condition ensuring a unique root on the positive real axis is $Bb-A > 0$, which is \eqref{eq:bvcond}. Such a condition implies the existence of a unique flexus on $\gamma>0$: this follows from the fact that 
\begin{equation}
\lim_{\gamma \rightarrow +\infty} F(\gamma) = B+C,\notag
\end{equation}
but the claim can be also verified by applying Descartes' rule on $\frac{d^2F(\gamma)}{d^2\gamma}$. \\Finally, the optimal $\gamma$ is obtained by solving \eqref{eq:firstDerBV} and has the following expression:
\begin{equation}
    \hat{\gamma}_{(a)} = \frac{b(A + \sqrt{ABb})}{Bb-A}.\notag
\end{equation}\\
(b) We proceed along the lines of the preceding argument, but considering the approximation error bound as in \eqref{eq:approxerror2}. Considering the following coefficients:
\begin{equation}
    A \text{ as in \eqref{eq:Abbv},} \qquad D = \sum_{i=1}^E \frac{\|\bar{\alpha}^{\pi}\|_{\infty}}{\lambda_i},
\end{equation}
the claim follows by proving that the function $F(\gamma) = \frac{A}{\gamma} + D\gamma$ has a unique minimum for $\gamma >0$. This is shown by studying the first and second derivatives, and using the fact that both $A$ and $D$ are positive. \\
The resulting optimal $\gamma$ always exists and takes the following value:
\begin{equation}
    \hat{\gamma}_{(b)} = \sqrt{\frac{A}{D}}.\notag
\end{equation}

\subsection{Statement of benchmark Theorems in Section \ref{sec:discussion}}\label{sec:statementsdisc}
For ease of referencing, we report the statements of \cite[Theorem 5]{Smale2007LearningTE} and \cite[Proposition 20]{Lin2017distributed} used in Sections \ref{sec:discSZ} and \ref{sec:discLGZ}, respectively.
\begin{theorem}[\cite{Smale2007LearningTE}, Theorem 5]
Let $\rho$ satisfy $|y| \leq M$ almost surely. Then for any $0 < \delta < 1$ 
\begin{equation}
  \|f_z - f_{\cH}\|_{\cL} \leq \frac{12 C_{\cK}M \log(4/\delta)}{\sqrt{N\gamma}}
    \label{eq:SZbound}
\end{equation}
provided that
\begin{equation}
\gamma \geq \frac{8C_{\cK}^2\log(4/\delta)}{\sqrt{N}}.
    \label{eq:SZgamma}
\end{equation}
\hfill$\square$
\end{theorem}

\begin{proposition}[\cite{Lin2017distributed}, Proposition 20] 
Assume $\mathbb{E}[y^2] < \infty$ and that $\sigma_{\rho}^2 \in \mathscr{L}^p_{\cX}$ for some $1 \leq p \leq \infty$. Moreover, define the effective dimension $\mathcal{N}(\gamma) = \textup{Tr}((L_{\cK} + \gamma I)^{-1}L_{\cK})$~\cite{zhang2005effdim}. Then, 
\begin{align}\label{eq:LGZexpval}
    \mathbb{E}_{\cZ} [\|f_z - f_{\cH}\|_{\cL}] \leq& (2 + 56 C_{\cK}^4 + 57C_{\cK}^2) \Big(1 + \frac{1}{(N\gamma)^2} + \frac{\mathcal{N}(\gamma)}{N\gamma}\Big) \notag \\
    &\times \Bigg\lbrace C_{\cK}^{\frac{1}{p}}\sqrt{\|\sigma_{\rho}\|_p}\Big(\frac{\mathcal{N}(\gamma)}{N}\Big)^{\frac{1}{2}(1-\frac{1}{p})}\Big(\frac{1}{N\gamma}\Big)^{\frac{1}{2p}} + \notag  + C_{\cK}\frac{\|f_{\cH} - f_{\rho}\|_{\cL}}{\sqrt{N\gamma}} \Bigg\rbrace.
\end{align}
\hfill$\square$
\end{proposition}

\subsection{Comparison with sample error bound in \cite{Wang2011}}\label{sec:discGZ}
We now perform another Monte Carlo study in another set-up where noises have unbounded support. We compare bound~\eqref{eq:uniformbound} with the following:
\begin{theorem}[\cite{Wang2011}, Theorem 1]
Assume there exist constants $\tilde{M}>0$ and $C>0$ such that the so-called moment hypothesis holds:
\begin{equation}
    \int_{\cY}|y|^{\ell}d\rho(y|x) \leq C\ell!\tilde{M}^{\ell} \qquad \forall \ell \in \mathbb{N},\, x\in\cX. \label{eq:momenthypothesis}
\end{equation}
Furthermore, assume that there exists some $0<\beta\leq 1$ and a constant $C_{\beta} > 0$ such that
\begin{equation}
    \|f_{\cH} - f_{\rho}\|_{\cL}^2 + \gamma\|f_{\cH}\|_{\cH}^2 \leq C_{\beta}\gamma^{\beta}. \label{eq:Cbeta}
\end{equation}
Then, if the kernel $\cK$ is infinitely differentiable on $\cX\times \cX$, then for any $0<\varepsilon<1$ and $0<\delta<1$, with confidence $1-\delta$ we have, by taking $\gamma = N^{\varepsilon - 1}$,
\begin{equation}
    \|f_z - f_{\cH}\|_{\cL}^2 \leq \tilde{C}_{\varepsilon}N^{\varepsilon-1}\log(4/\delta)^{\frac{4}{\varepsilon} + 2}.\label{eq:GZbound}
\end{equation}
\hfill$\square$
\label{theorem:GZ}
\end{theorem}

Before presenting the details of the numerical experiment, we derive the expressions for $C$, $\tilde{M}$, $C_{\beta}$ and $\tilde{C}_{\varepsilon}$. Theorem \ref{theorem:GZ} is a corollary of the following result:
\begin{theorem}[\cite{Wang2011}, Theorem 2]
Assume the moment hypothesis \eqref{eq:momenthypothesis} with constants $C$ and $\tilde{M}$ holds; moreover, let condition \eqref{eq:Cbeta} with $0<\beta \leq 1$ and constant $C_{\beta}$ be valid. Define $\mathscr{N}(\mathcal{B}_1,\eta)$ the minimum number of disks with radius $\eta$ that cover the balls $\mathcal{B}_1 = \{f \in \cH \: : \: \|f\|_{\cH} \leq 1 \}$, and assume that $\cH$ has polynomial complexity exponent $s>0$, i.e., 
\begin{equation}
    \log\mathscr{N}(\mathcal{B}_1,\eta) \leq C_0 \Big(\frac{1}{\eta}\Big)^s.\label{eq:polycomplexity}
\end{equation}
If $0 < \varepsilon < \frac{\beta}{1+s}$, then by taking $\gamma = N^{\frac{\varepsilon}{\beta} - \frac{1}{s+1}}$, for any $0<\delta<1$, with confidence $1-\delta$ we have
\begin{equation}
    \|f_z - f_{\rho}\|_{\cL}^2 \leq \tilde{C}_{\varepsilon}N^{\varepsilon - \frac{\beta}{s+1}}\Big(\log\Big(\frac{4}{\delta}\Big)\Big)^{\frac{\beta(1+\beta)}{(s+1)\varepsilon}+2}.
\end{equation}
The constant $\tilde{C}_{\varepsilon}$ can be computed as follows:
\begin{equation}
    \tilde{C}_{\varepsilon} = \frac{C_5}{\varepsilon^2}C_2^{\frac{4\beta}{\varepsilon(s+1)}}\Big(1 + \log\Big(1 + \frac{2}{\varepsilon(s+1)}\Big) \Big)^{\frac{\beta(1+\beta)}{\varepsilon(s+1)}+2}, \notag
\end{equation}
where we have
\begin{equation}
    \begin{cases}
    C_5 &= 38C_{\beta} + 2(C_1 + 32^2(C+1)^2)C_4^2(2/(s+1))^2 + 480(C_{\cK}+1)^2C_{\beta}\\
    C_4 &= \tilde{M}(2C_{\cK}(C + (1 + 2\sqrt{2C}) +1)) + C_3\\
    C_3 &= \sqrt{38C_{\beta}} + (C_{\cK} + 1)\sqrt{480 C_{\beta}} + \tilde{M}\\
    C_2 &= \sqrt{2[C_1 + 32^2(C+1)^2]}\\
    C_1 &= 6C_{\cK} + 6C + 8(1 + \sqrt{2C})/\tilde{M} + 520 (C_{\cK} + C + 2(C+1))^2(C_0+1).
    \end{cases}
\end{equation}
\hfill $\square$
\label{theorem:completeGZ}
\end{theorem}

The choice in \cite{Wang2011} to obtain the statement of Theorem \ref{theorem:GZ} from the result above presented is to assume that $f_{\rho}$ belongs to the hypothesis space $\cH$, and setting $\beta=1$ and $s = \frac{\varepsilon}{1 - \varepsilon}$ with $0 < \varepsilon < 1/2$ (and then rescaling $2\varepsilon$ to $\varepsilon$). Our task is now to obtain explicit values for the constants $C$, $\tilde{M}$, $C_{\beta}$ and $C_0$. The bound for $C_{\cK}$ is given in \eqref{eq:boundCk}.

\paragraph*{Constants for moment hypothesis \eqref{eq:momenthypothesis} ($C$ and $\tilde{M}$)}
Following Example 1 in \cite{Wang2011}, we obtain that the moment hypothesis \eqref{eq:momenthypothesis} is satisfied for $C=4$ and $\tilde{M} = \max\{\sqrt{B_0},B_{\infty}\}$, where $B_0 = B_{\sigma}^2$ (see Assumption \ref{assume:bounds}) and $B_{\infty} \geq \|f_{\rho}\|_{\infty}$. Note that $B_{\infty} < \infty$ because $f_{\rho}$ belongs to $\cH$.

\paragraph*{Constant entering bound \eqref{eq:Cbeta} ($C_{\beta}$)}
We evaluate the cost \eqref{eq:probdatafree} at $f_{\cH} = (I + \gamma\io^{-1})^{-1}f_{\rho}$. We then obtain
\begin{align*}
    \|f - f_{\rho}\|_{\cL}^2 + \gamma\|\io^{-1/2}f\|_{\cL}^2 &= \|[(I + \gamma\io^{-1})^{-1} - I]f_{\rho}\|_{\cL}^2 + \gamma\|\io^{-1/2}(I + \gamma\io^{-1})^{-1}f_{\rho}\|_{\cL}^2\\
    &=\sum_{i=1}^E \Big[\Big(\frac{\lambda_i}{\lambda_i + \gamma} - 1\Big)^2 + \gamma\Big(\frac{\sqrt{\lambda_i}}{\lambda_i + \gamma} \Big)^2 \Big](\bar{\alpha}_i^{\pi})^2\\
    &= \gamma\sum_{i=1}^E \frac{1}{\lambda_i + \gamma}(\bar{\alpha}_i^{\pi})^2 = \gamma\sum_{i=1}^E \Big( \frac{\lambda_i^{\beta}}{\lambda_i + \gamma}\Big)\lambda_i^{-\beta}(\bar{\alpha}_i^{\pi})^2
\end{align*}
for any $\beta$ in $(0,1]$. We can bound such an expression as follows (see Theorem 3, Chapter II.2 in \cite{Cucker02} for the whole derivation):
\begin{equation}
   \gamma\sum_{i=1}^E \Big( \frac{\lambda_i^{\beta}}{\lambda_i + \gamma}\Big)\lambda_i^{-\beta}(\alpha_i^{\pi})^2 \leq \gamma \Big(\sup_{\tau} \frac{\tau^{\beta}}{\tau + \gamma} \Big)\|\io^{-\beta/2}f_{\rho}\|_{\cL}^2 \leq \gamma^{\beta}\|\io^{-\beta/2}f_{\rho}\|_{\cL}^2.\notag
\end{equation}
Therefore, the constant $C_{\beta}$ is obtained by bounding $\|\io^{-\beta/2}f_{\rho}\|_{\cL}^2$. In the case $\beta=1$, which is of our interest, a solution could be $C_{\beta=1} \leq B_f^2/\lmin$.

\paragraph*{Constant for polynomial complexity \eqref{eq:polycomplexity} ($C_0$)}
Since $\cK$ in \eqref{eq:ourkernel} is infinitely differentiable, condition \eqref{eq:polycomplexity} is known to hold for any $s > 0$; however, since we are considering $s = \frac{\varepsilon}{1 - \varepsilon}$ with $0 < \varepsilon < 1/2$, we are restricting our attention to $0<s<1$.\\ 
Since $\cH$ is finite-dimensional with dimension $E$, Theorem 5.3 \cite{cuckerzhou2007} gives that $\mathscr{N}(\mathcal{B}_1,\eta) \leq (1 + 2/\eta)^E$. Taking logarithms on both sides, we obtain
\begin{equation*}
    \log\mathscr{N}(\mathcal{B}_1,\eta) \leq E\log\Big(1 + \frac{2}{\eta}\Big) \leq 2EG(s)\Big(\frac{1}{\eta}\Big)^s.
\end{equation*}
We now need to find $G(s)$. Its expression is summarized in the following result.
\begin{lemma}
Consider $x>0$. In order to have $\log(1+x) < G(s)x^s$ for $0<s<1$, it has to hold that $$G(s) > \Big(\frac{1-s}{s}\Big)^{1-s}.$$
\end{lemma}
\textit{Proof. } 
Consider the function $g(x) = G(s)x^s - \log(1+x)$. We want it to be always positive for $x>0$. Since $g(0) = 0$, this amounts to imposing that $\frac{dg(x)}{dx}$ is always positive. 
So, studying the first derivative, we have
\begin{equation}
\begin{aligned}
    &\frac{dg(x)}{dx} = \frac{G(s)s(x+1)-x^{1-s}}{x^{1-s}(1+x)} > 0 \longrightarrow \frac{x}{(1+x)^{\frac{1}{1-s}}} < \Big(G(s) s\Big)^{\frac{1}{1-s}}.
    \end{aligned}
    \notag
\end{equation}
The claim follows by maximising the term on the left-hand side.
\hfill $\blacksquare$

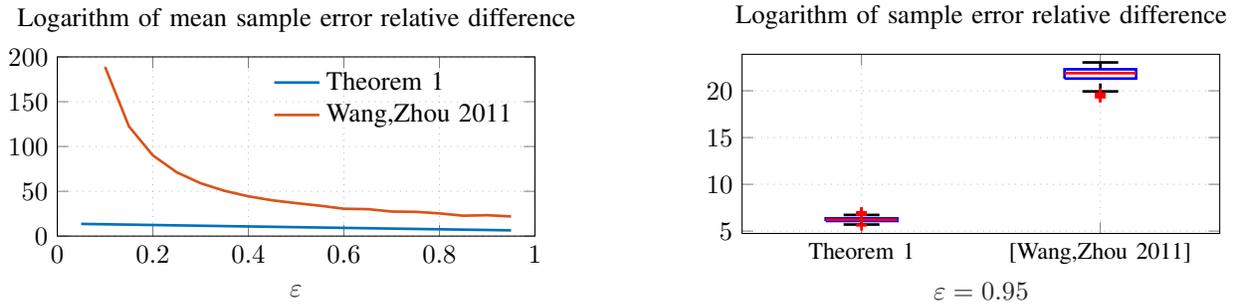
\begin{figure}[h!]
    \begin{minipage}{0.5\textwidth}
    \centering
    \input{compWZ2.tex}
    \end{minipage}
    \begin{minipage}{0.5\textwidth}
    \centering
    \input{compWZ1.tex}
    \end{minipage}
    \caption{Statistics of the difference between bound and true sample error, normalized by true sample error, in logarithmic scale. Top panel: behaviour of the mean value over the Monte Carlo iterations as a function of $\varepsilon$. Bottom panel: boxplots over the Monte Carlo runs for the case $\varepsilon=0.95$. The values (in logarithmic scale) returned by Theorem \ref{prop:uniformbound} are $6.238 \pm 4.694$, while the ones obtained using Theorem \ref{theorem:GZ} are $21.961 \pm 21.421$.}
    \label{fig:WZ}
\end{figure}

We are now ready to perform the numerical test. We randomize both on the regression function and on the number of data-points. The first is drawn in the same way as in Section~\ref{sec:discSZ}, while data-set cardinalities $N$ 
are drawn from the set $\{300, 315,...,6990\}$. As in Section~\ref{sec:discSZ}, the SNR is set to 150, but noises are now distributed as Gaussian. The hyperparameters $\{\lambda_i\}_{i=1}^{20}$ 
ruling both the sampling of $\bar{\alpha}_q$ and the hypothesis space $\cH$ are all set to 10. The numerical values of the sample error bounds are computed with $\gamma$ set as in the statement of Theorem \ref{theorem:GZ}. For each Monte Carlo run, we evaluate the bounds with $\varepsilon$ in the grid $\{0.05,0.1,...,0.95\}.$ To compare the numerical values of the bounds, the score we consider is the difference between bound and true sample error, divided by the true sample error. The results are displayed in Figure \ref{fig:WZ}. From~\eqref{eq:GZbound} it is clear that the convergence rate in the number of data $N$ is better compared with that presented in~\eqref{eq:uniformbound}; however, the numerical values returned by the choice of $\gamma$ are extremely conservative.

\subsection{Additional results for the tests in Section \ref{subsec:bv}}\label{sec:addendaBV}
We recall that the numerical values chosen for the Monte Carlo test were the following:
\begin{itemize}
    \item $\cX = [-5\times 10^4,\,5\times 10^4]$;
    \item $\delta = 0.5$;
    \item Regression function $f_{\rho} = \sum_{q}\bar{\alpha}_q\bar{\varphi}_q$ characterized by 30 sine/cosine pairs $\{\bar{\varphi}\}_{q=1}^{30}$, where each $q$ is randomly selected without repetitions from the set $\{1,...,100\}$;
    \item all components of $\bar{\alpha} \in \mathbb{R}^{30}$ are independent samples from a Gaussian distribution with zero mean and variance $\lambda=1$;
    \item in the hypothesis space, $\lambda_i = \lambda$ for all $i=1,...,E$;
    \item $Q$ is selected as a random subset with cardinality 10 from the set of frequencies characterising the regression function;
    \item SNR = 50;
    \item at each iteration, we draw a random regression function and select the basis functions of $\cH$;
    \item $N = 2500.$
\end{itemize}

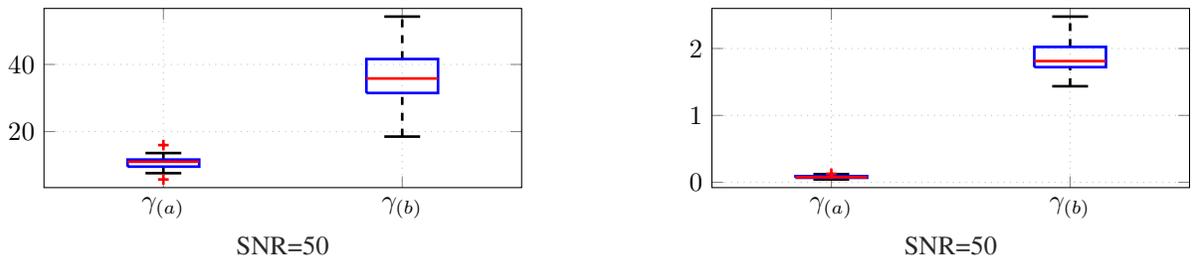
\begin{figure}[H]
   \begin{minipage}{0.5\textwidth}
   \centering
   \input{BV_sample.tex}
   \end{minipage}
   \begin{minipage}{0.5\textwidth}
   \centering
   \input{BV_approx.tex}
   \end{minipage}
    \caption{Boxplots of the normalized, relative differences between bounds and true values. Left panel: results for the sample error; right panel: results for the approximation errors. These complement the results presented in Table \ref{tab:bv}, Section \ref{subsec:bv}.}
    \label{fig:BVs}
\end{figure}

We further display the results reported in Table \ref{tab:bv} in Figure \ref{fig:BVs}. In each boxplot, we consider the difference between the bound and the true error, and we normalize it by the latter.

\subsection{The benefits of regularization: the case of additive noise model} \label{subsec:regularization}
\begin{figure}[h]
    \centering
    \input{analysis_reg.tex}
    \caption{Relative discrepancy $100\%(\bar{\mathcal{M}}(\hat{\gamma}) - \bar{\mathcal{M}}(\gamma^*))/\bar{\mathcal{M}}(\gamma^*)$ for the regularized and non-regularized cases. } 
    \label{fig:regularization}
\end{figure}
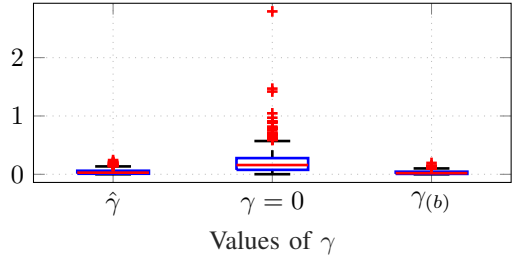

We now discuss the performance of the estimation scheme proposed in \eqref{eq:fz} with respect to the one that would have been obtained without regularization. We carry out the analysis assuming an additive noise model, and regarding the regression function as the "true" function to be estimated from data. In this setting, the measurements model for each $t=1,...,N$ is $y_t = \phi^{\top}(x_t)\bar{\alpha}^{\pi} + r(x_t) + e_t$, where the first term is given by the projection of the regression function on the subspace spanned by the basis functions $\{\varphi_i(\cdot) \}_i\}_{i=1}^E$ defined in \eqref{eq:basisRKHS} and entering $\phi^{\top}(\cdot)$; the second collects the contribution of frequencies that have not been included in the subspace and is assumed to be bounded; the last one is the additive noise, which is assumed to be i.i.d.~with zero mean and known variance $\sigma^2$. The overall estimation problem to be solved reads as follows:
\begin{equation}
    \hat{\alpha} = \arg\min_{\alpha \in \mathbb{R}^E} \|Y - \Phi \alpha\|^2 + N\sigma^2 \alpha^{\top}P^{-1}\alpha, \label{eq:alphahat}
    \end{equation}
where $\Phi \in \mathbb{R}^{N\times E}$ stacks all $\{\phi^{\top}(x_t)\}_{t=1}^N$, $Y$ is defined as the one entering \eqref{eq:solSampOp}, and $P^{-1}$ is a regularization matrix. The solution \eqref{eq:alphahat} reads as 
\begin{equation}
    \hat{\alpha} = 
    (\Phi^{\top}\Phi + N\sigma^2 P^{-1})^{-1}\Phi^{\top}Y. \label{eq:hatalpha}
\end{equation}
Our goal is to discuss the performance of such an estimator when considering $P^{-1} = \frac{\gamma \Sigma_{\alpha}^{-1}}{\sigma^2}$, which directly relates to \eqref{eq:fz}, as a function of $\gamma$. To do so, we consider $\mathcal{M}(\gamma)=\mathbb{E}_e[\|\hat{\alpha} - \bar{\alpha}^{\pi}\|^2]$, where $\mathbb{E}_e[\cdot]$ denotes expectation with respect to the noise distribution, as a performance score. Carrying out the computations with the particular choice of the regularizer above introduced, and defining $r = [r(x_1), ..., r(x_N)]^{\top}$ and $\tilde{\Sigma}_{\alpha} = \Sigma_{\alpha}/N$, it turns out that
\begin{align}
    \mathcal{M}(\gamma&) = \text{Tr}\Bigg((\gamma\tilde{\Sigma}_{\alpha}^{-1} + \Phi^{\top}\Phi)^{-1}\Big[ \sigma^2\Phi^{\top}\Phi + \gamma^2\tilde{\Sigma}_{\alpha}^{-1}\bar{\alpha}^{\pi}(\bar{\alpha}^{\pi})^{\top}\tilde{\Sigma}_{\alpha}^{-1} + \Phi^{\top}rr^{\top}\Phi + 2\gamma\tilde{\Sigma}_{\alpha}^{-1}\bar{\alpha}^{\pi}r^{\top}\Phi\Big](\gamma\tilde{\Sigma}_{\alpha}^{-1} + \Phi^{\top}\Phi)^{-1} \Bigg).
    \label{eq:MSEalpha}
\end{align}
Note that $\mathcal{M}(0) = \text{Tr}((\Phi^{\top}\Phi)^{-1}[\Phi^{\top}rr^{\top}\Phi + \sigma^2\Phi^{\top}\Phi](\Phi^{\top}\Phi)^{-1})$. 
Denoting with $R$ the expression in square brackets in \eqref{eq:MSEalpha}, we obtain
\begin{align}
    \frac{d\mathcal{M}(\gamma)}{d\gamma}=2\text{Tr}\Bigg((\gamma\tilde{\Sigma}_{\alpha}^{-1} + \Phi^{\top}\Phi)^{-1}\Big[\gamma\tilde{\Sigma}_{\alpha}^{-1}\bar{\alpha}^{\pi}(\bar{\alpha}^{\pi})^{\top}\tilde{\Sigma}_{\alpha}^{-1} + \tilde{\Sigma}_{\alpha}^{-1}\bar{\alpha}^{\pi}r^{\top}\Phi  - \tilde{\Sigma}_{\alpha}^{-1}(\gamma\bar{\Sigma}_{\alpha}^{-1} + \Phi^{\top}\Phi)^{-1}R \Big] (\gamma\tilde{\Sigma}_{\alpha}^{-1} + \Phi^{\top}\Phi)^{-1}\Bigg).
    \notag
\end{align}
Studying its limit as $\gamma \rightarrow 0^+$, one has
\begin{align}
2\text{Tr}\Bigg((\Phi^{\top}\Phi)^{-1}\tilde{\Sigma}_{\alpha}^{-1}\Big[\bar{\alpha}^{\pi}r^{\top}\Phi - (\Phi^{\top}\Phi)^{-1}(\Phi^{\top}rr^{\top}\Phi + \sigma^2\Phi^{\top}\Phi) \Big] (\Phi^{\top}\Phi)^{-1}\Bigg).\label{eq:dMat0}
\end{align}
The only case that is easy to study occurs when there is no residual term (i.e., the regression function belongs to the hypothesis space, so that $r(x) = 0$ for all $x$). In that scenario, 
the expression above is clearly negative, and this proves the fact that $\mathcal{M}(\gamma) < \mathcal{M}(0)$ at least in some small neighbourhood of the origin: see also Proposition 2 in \cite{MU2018381}.\\ 
To numerically test the impact of regularization, we perform a Monte Carlo study of 500 trials. We consider an input domain $\cX = [-25,25]$, and a regression function characterized by 50 basis functions randomly selected among the first 80 (ordered with increasing $q$ in \eqref{eq:basis}), linearly combined by a vector drawn from a Gaussian distribution with zero mean, i.i.d.~components and variance $\lambda=10$. The latter hyper-parameter also enters the definition of $\cH$; the set of frequencies $Q$ defining it is a sample of random dimension $E/2$ in $\{5,6,...,50\}$, selected among the ones defining $f_{\rho}$. Thus, in this testing situation the approximation error ruled by the residual $r(\cdot)$ is different from 0. We consider an SNR of 100 yielded by a Gaussian, zero-mean, i.i.d.~noise. On each run, the number of data $N$ is randomly selected in the interval $\{5,...,E/2\}$. For the regularization parameter, we both use $\gamma^{(b)}$ as in Proposition \ref{prop:tradeoff}(b) and a value $\hat{\gamma}$ estimated via marginal likelihood optimization performed with a Gibbs sampling scheme leveraging the Bayesian interpretation of the problem in \eqref{eq:alphahat} (\cite[Chapter 1]{gilks1995markov}; see also \cite{PILLONETTO2015106} for a thorough discussion on the robustness of the marginal likelihood hyper-parameter estimation. More details are reported below). We then study the values of the overall error $\bar{\mathcal{M}}(\gamma) = \|f_{\rho} - f_z\|_{\cL}$, where $f_z(\cdot) = \phi^{\top}\hat{\alpha}$ with $\hat\alpha$ as in \eqref{eq:hatalpha}, attained by $\gamma^{(b)}$, $\gamma=0$ and $\hat{\gamma}$, compared to the oracle value corresponding to $\gamma=\gamma^*$ computed by grid search. The results are summarized in Figure \ref{fig:regularization}. We note that the two regularized estimators yield comparable results, meaning that $\gamma^{(b)}$ is a good estimator of $\hat{\gamma}$ (and is faster to be computed); and that both always outperform the case with no regularization involved. Specifically, we obtain that the relative discrepancy between, e.g., the case with $\gamma=0$ and $\hat{\gamma}$ has a median value of $24.99\%$, with minimum and maximum values equal to $0.64\%$ and $66.97\%$, respectively.\\

We conclude by providing the details for the selection rule for hyper-parameter $\gamma$ used as baseline in the Monte Carlo test.  Let us consider the estimation problem in \eqref{eq:alphahat}, i.e.,
\begin{align}
    \hat{\alpha} = \arg\min_{\alpha \in \mathbb{R}^E} \|Y - \Phi \alpha\|^2 + N\sigma^2 \alpha^{\top}P^{-1}\alpha  
     =(\Phi^{\top}\Phi + N\sigma^2 P^{-1})^{-1}\Phi^{\top}Y,  \label{eq:alphahatAPP} 
    \end{align}
where $\Phi \in \mathbb{R}^{N\times E}$ is the matrix that stacks all $\{\phi^{\top}(x_t)\}_{t=1}^N$ (i.e., the row vectors containing the basis functions $\{\varphi_i(\cdot)\}_{i=1}^E$ as defined in \eqref{eq:basisRKHS}), $P^{-1}$ is a regularization matrix, and $Y$ is the output measurements vector $[y_1,\,...,y_N]^{\top}$.\\ The objective in \eqref{eq:alphahatAPP} admits a stochastic interpretation. Consider a measurements model $Y = \Phi\alpha + E$, where $E$ is an $N$-dimensional Gaussian vector with zero mean and known covariance $\sigma^2 I_N$, and $\alpha \in \mathbb{R}^E$ is unknown. Taking the Bayesian viewpoint, $\alpha$ is modelled as a random vector; specifically, assume it is Gaussian, with zero mean and covariance $P/N$. In this set-up, the objective in \eqref{eq:alphahatAPP} is (apart from constants not depending on $\alpha$) the negative logarithm of the posterior probability $\alpha|Y$: that is, $\hat{\alpha}$ is computed as the Maximum a Posteriori (MAP) estimate, which corresponds to the minimum variance linear estimate in the Gaussian case we are considering. In fact, the solution to \eqref{eq:alphahatAPP} is indeed the expression of the posterior mean. \\
Within this framework, assume that (some of the) hyper-parameters entering $P^{-1}$, or $\sigma^2$, are unknown, and collect them in a vector
$\eta$. A possible strategy consists in estimating them from data, leveraging the so-called empirical Bayes approach~\cite{maritz2018empirical}: in particular, the estimate for $\eta$ is computed by maximising the evidence $Y|\eta$ (or, more conveniently, minimizing its negative logarithm), which is equivalent to the joint distribution $Y, \alpha|\eta$ where the dependence from $\alpha$ is integrated out. In the Gaussian case, this reads as
\begin{equation}
    \hat{\eta} = \arg\min_{\eta}\: Y^{\top}\Sigma_Y(\eta)^{-1}Y + \log\det\Sigma_Y(\eta),\label{eq:marglikh}
\end{equation}
where $\Sigma_Y = \Phi P\Phi^{\top}/N + \sigma^2I_N$. This problem is non-convex, and deterministic optimization routines might return unreliable results due to their sensitivity to initial conditions. A way to overcome this issue consists in resorting to Markov Chain Monte Carlo (MCMC). Specifically, the idea is to run a (single-component) Metropolis-Hastings algorithm to construct a Markov chain whose invariant distribution is (proportional to) the marginal likelihood of interest: this is a mechanism to draw samples from such a distribution, and solve \eqref{eq:marglikh} in sample-based form. For an introduction to MCMC we refer to \cite{gilks1995markov}.\\

Let us now relate \eqref{eq:alphahatAPP} to the original function estimation problem stated in \eqref{eq:fz} to provide an expression for $P^{-1}$, and detail the MCMC-based procedure for marginal likelihood optimization. The function is estimated as
\begin{equation}
    f_z = \arg\min_{f \in \cH} \frac{1}{N}\sum_{t=1}^N (y_t - f(x_t))^2 + \gamma \|f\|_{\cH}^2. \notag
\end{equation}
According to the definitions \eqref{eq:RKHS} and \eqref{eq:ourkernel} specifying the RKHS $\cH$, one can write $f_z(\cdot) = \phi^{\top}(\cdot)\hat{\alpha}$. In this way, estimating $f_z(\cdot)$ translates into solving \eqref{eq:alphahatAPP} for the particular choice $P^{-1} = \frac{\gamma\Sigma_{\alpha}^{-1}}{\sigma^2}$, where $\Sigma_{\alpha}$ is given by the choice of the hypothesis space (see, e.g., \eqref{eq:ourkernel}), and $\gamma$ is the positive scalar to be tuned.\\
Let us now detail the MCMC procedure for the scenario above specified. First, denoting with $p(\cdot)$ the probability density function of interest, we note that 
\begin{equation}
    p(\alpha,\gamma|Y) = p(\alpha|\gamma,Y)p(\gamma|Y) \propto p(\alpha | \gamma, Y)p(Y|\gamma), \notag
\end{equation}
which tells us that drawing samples from $\alpha,\gamma|Y$ is a suitable way to explore the marginal likelihood. Therefore, we set up a Gibbs sampler (a particular case of single-component Metropolis-Hastings algorithm) whose invariant distribution has density $p(\alpha,\gamma|Y) = \pi(\alpha,\gamma)$. In this particular scenario, samples from the full conditionals are easy to be computed: in fact, denoting with $\tilde{\Sigma}_{\alpha} = \Sigma_{\alpha}/N$ and with $\Gamma(\mathfrak{a},\mathfrak{b})$ a Gamma distribution with mean $\mathfrak{a}/\mathfrak{b}$, using conjugate distributions properties we obtain
\begin{align}
    &\pi(\alpha|\gamma) \longleftarrow \alpha|\gamma,Y \sim \mathcal{N}\Big( (\tilde{\Sigma}_{\alpha}^{-1} + \Phi^{\top}\Phi/\gamma)^{-1}\Phi^{\top}Y/\gamma,\: (\tilde{\Sigma}_{\alpha}^{-1} + \Phi^{\top}\Phi/\gamma)^{-1} \Big)\\
    &\pi(\gamma|\alpha) \longleftarrow \gamma|\alpha,Y \sim \Gamma\Big(\frac{N}{2},\: \frac{\|Y - \Phi\alpha\|^2}{2}\Big).
\end{align}
Finally, the last $N_g < N_G$ samples obtained from the Gibbs sampler can be used to compute $\hat{\gamma}$ maximising the marginal likelihood.

\end{document}

%% file: SZandLGZcomparison.tex
%
%
\begin{tikzpicture}

\begin{axis}[%
scale = 0.25,
width=12.0in,
height=3.754in,
at={(1.136in,0.642in)},
scale only axis,
xmin=0.5,
xmax=3.5,
xtick={1,2,3},
xticklabels={{Theorem 1},{[Smale, et al. 2007]},{[Lin, et al. 2017]}},
xticklabel style={font = \small},
ymin=1.76142703778331,
ymax=17.2356183735049,
axis background/.style={fill=white},
title style={font=\normalsize},
title={Logarithm of sample error relative difference},
xmajorgrids,
ymajorgrids,
grid style={dotted}
]
\addplot [color=black, dashed, line width=1.0pt, forget plot]
  table[row sep=crcr]{%
1	3.16012883658933\\
1	3.51281537089617\\
};
\addplot [color=black, dashed, line width=1.0pt, forget plot]
  table[row sep=crcr]{%
2	6.22099058841271\\
2	6.65320054769995\\
};
\addplot [color=black, dashed, line width=1.0pt, forget plot]
  table[row sep=crcr]{%
3	15.5450881911593\\
3	16.532246040063\\
};
\addplot [color=black, dashed, line width=1.0pt, forget plot]
  table[row sep=crcr]{%
1	2.60635131434642\\
1	2.92484417245146\\
};
\addplot [color=black, dashed, line width=1.0pt, forget plot]
  table[row sep=crcr]{%
2	5.46805070747596\\
2	5.90850768476421\\
};
\addplot [color=black, dashed, line width=1.0pt, forget plot]
  table[row sep=crcr]{%
3	10.0851272453475\\
3	13.3549338306358\\
};
\addplot [color=black, line width=1.0pt, forget plot]
  table[row sep=crcr]{%
0.8875	3.51281537089617\\
1.1125	3.51281537089617\\
};
\addplot [color=black, line width=1.0pt, forget plot]
  table[row sep=crcr]{%
1.8875	6.65320054769995\\
2.1125	6.65320054769995\\
};
\addplot [color=black, line width=1.0pt, forget plot]
  table[row sep=crcr]{%
2.8875	16.532246040063\\
3.1125	16.532246040063\\
};
\addplot [color=black, line width=1.0pt, forget plot]
  table[row sep=crcr]{%
0.8875	2.60635131434642\\
1.1125	2.60635131434642\\
};
\addplot [color=black, line width=1.0pt, forget plot]
  table[row sep=crcr]{%
1.8875	5.46805070747596\\
2.1125	5.46805070747596\\
};
\addplot [color=black, line width=1.0pt, forget plot]
  table[row sep=crcr]{%
2.8875	10.0851272453475\\
3.1125	10.0851272453475\\
};
\addplot [color=blue, line width=1.0pt, forget plot]
  table[row sep=crcr]{%
0.775	2.92484417245146\\
0.775	3.16012883658933\\
1.225	3.16012883658933\\
1.225	2.92484417245146\\
0.775	2.92484417245146\\
};
\addplot [color=blue, line width=1.0pt, forget plot]
  table[row sep=crcr]{%
1.775	5.90850768476421\\
1.775	6.22099058841271\\
2.225	6.22099058841271\\
2.225	5.90850768476421\\
1.775	5.90850768476421\\
};
\addplot [color=blue, line width=1.0pt, forget plot]
  table[row sep=crcr]{%
2.775	13.3549338306358\\
2.775	15.5450881911593\\
3.225	15.5450881911593\\
3.225	13.3549338306358\\
2.775	13.3549338306358\\
};
\addplot [color=red, line width=1.0pt, forget plot]
  table[row sep=crcr]{%
0.775	3.044894582244\\
1.225	3.044894582244\\
};
\addplot [color=red, line width=1.0pt, forget plot]
  table[row sep=crcr]{%
1.775	6.05369012144661\\
2.225	6.05369012144661\\
};
\addplot [color=red, line width=1.0pt, forget plot]
  table[row sep=crcr]{%
2.775	14.6553746552207\\
3.225	14.6553746552207\\
};
\addplot [color=black, line width=1.0pt, only marks, mark=+, mark options={solid, draw=red}, forget plot]
  table[row sep=crcr]{%
1	2.4647993712252\\
1	2.55867223125809\\
1	2.56352950132624\\
1	2.56577139304885\\
1	2.56809528191386\\
1	3.52689370021957\\
1	3.52960672133012\\
1	3.54840167801601\\
1	3.57983202507436\\
1	3.61554230579912\\
};
\addplot [color=black, line width=1.0pt, only marks, mark=+, mark options={solid, draw=red}, forget plot]
  table[row sep=crcr]{%
2	6.68974796226307\\
};
\addplot [color=black, line width=1.0pt, only marks, mark=+, mark options={solid, draw=red}, forget plot]
  table[row sep=crcr]{%
3	5.53890325987907\\
3	6.45497170029377\\
3	7.41770186053894\\
3	7.70687928475331\\
3	7.71107214722567\\
3	7.84333996370317\\
3	8.23059780902931\\
3	8.37879017060551\\
3	8.49322121812499\\
3	8.65549958774682\\
3	8.76274821063522\\
3	8.89738421348669\\
3	8.99356761672983\\
3	9.28638148982074\\
3	9.29211507629291\\
3	9.41160085538785\\
3	9.43595393949429\\
3	9.56211643950449\\
3	9.61958852164132\\
3	9.79843488501958\\
3	9.92150604216283\\
3	9.985753205445\\
};
\end{axis}

\end{tikzpicture}%

%% file: compWZ2.tex
%
%
\definecolor{mycolor1}{rgb}{0.00000,0.44700,0.74100}%
\definecolor{mycolor2}{rgb}{0.85000,0.32500,0.09800}%
\begin{tikzpicture}

\begin{axis}[%
scale = 0.25,
width=10.0in,
height=3.754in,
at={(1.011in,0.8in)},
scale only axis,
unbounded coords=jump,
xmin=0,
xmax=1,
xlabel style={font=\color{white!15!black}},
xlabel={$\varepsilon$},
ymin=0,
ymax=200,
axis background/.style={fill=white},
title style={font=\normalsize},
title={Logarithm of mean sample error relative difference},
xmajorgrids,
ymajorgrids,
grid style={dotted},
legend style={legend cell align=left, align=left, fill=none, draw=none}
]
\addplot [color=mycolor1, line width=1.0pt]
  table[row sep=crcr]{%
0.05	13.6001145619612\\
0.1	13.1907477174236\\
0.15	12.7823615887422\\
0.2	12.374999497716\\
0.25	11.9687063041981\\
0.3	11.5635283463371\\
0.35	11.1595133595404\\
0.4	10.7567103721011\\
0.45	10.3551695755545\\
0.5	9.954942168045\\
0.55	9.55608016927894\\
0.6	9.15863620604218\\
0.65	8.76266326776271\\
0.7	8.36821443220815\\
0.75	7.97534256211912\\
0.8	7.58409997438287\\
0.85	7.19453808423096\\
0.9	6.80670702787605\\
0.95	6.42065526795548\\
};
\addlegendentry{Theorem 1}

\addplot [color=mycolor2, line width=1.0pt]
  table[row sep=crcr]{%
0.05	inf\\
0.1	188.893556467814\\
0.15	122.397616042032\\
0.2	90.1012689008088\\
0.25	71.2359470494001\\
0.3	59.0012990492556\\
0.35	50.5279685494882\\
0.4	44.4098890725814\\
0.45	39.9064762048287\\
0.5	36.7829368015798\\
0.55	33.923437995975\\
0.6	30.4933025061696\\
0.65	30.0674242177985\\
0.7	27.3925208847594\\
0.75	27.0745666605125\\
0.8	25.348643685941\\
0.85	22.7428369929948\\
0.9	23.3156891224709\\
0.95	21.9334806851029\\
};
\addlegendentry{Wang,Zhou 2011}

\end{axis}

\end{tikzpicture}%

%% file: compWZ1.tex
%
%
\begin{tikzpicture}

\begin{axis}[%
scale = 0.25,
width=10.0in,
height=3.754in,
at={(1.011in,0.8in)},
scale only axis,
xmin=0.5,
xmax=2.5,
xtick={1,2},
xticklabels={{Theorem 1},{[Wang,Zhou 2011]}},
xticklabel style={font = \small},
xlabel style={font=\color{white!15!black}},
xlabel={$\varepsilon = 0.95$},
ymin=4.74440689498115,
ymax=23.9017331555338,
axis background/.style={fill=white},
title style={font=\normalsize},
title={Logarithm of sample error relative difference},
xmajorgrids,
ymajorgrids,
grid style={dotted}
]
\addplot [color=black, dashed, line width=1.0pt, forget plot]
  table[row sep=crcr]{%
1	6.35633210567536\\
1	6.72587024119511\\
};
\addplot [color=black, dashed, line width=1.0pt, forget plot]
  table[row sep=crcr]{%
2	22.2978397043597\\
2	23.0309455982359\\
};
\addplot [color=black, dashed, line width=1.0pt, forget plot]
  table[row sep=crcr]{%
1	5.69836868498609\\
1	6.07916378383377\\
};
\addplot [color=black, dashed, line width=1.0pt, forget plot]
  table[row sep=crcr]{%
2	19.9355102267331\\
2	21.3080955154708\\
};
\addplot [color=black, line width=1.0pt, forget plot]
  table[row sep=crcr]{%
0.925	6.72587024119511\\
1.075	6.72587024119511\\
};
\addplot [color=black, line width=1.0pt, forget plot]
  table[row sep=crcr]{%
1.925	23.0309455982359\\
2.075	23.0309455982359\\
};
\addplot [color=black, line width=1.0pt, forget plot]
  table[row sep=crcr]{%
0.925	5.69836868498609\\
1.075	5.69836868498609\\
};
\addplot [color=black, line width=1.0pt, forget plot]
  table[row sep=crcr]{%
1.925	19.9355102267331\\
2.075	19.9355102267331\\
};
\addplot [color=blue, line width=1.0pt, forget plot]
  table[row sep=crcr]{%
0.85	6.07916378383377\\
0.85	6.35633210567536\\
1.15	6.35633210567536\\
1.15	6.07916378383377\\
0.85	6.07916378383377\\
};
\addplot [color=blue, line width=1.0pt, forget plot]
  table[row sep=crcr]{%
1.85	21.3080955154708\\
1.85	22.2978397043597\\
2.15	22.2978397043597\\
2.15	21.3080955154708\\
1.85	21.3080955154708\\
};
\addplot [color=red, line width=1.0pt, forget plot]
  table[row sep=crcr]{%
0.85	6.21949269173168\\
1.15	6.21949269173168\\
};
\addplot [color=red, line width=1.0pt, forget plot]
  table[row sep=crcr]{%
1.85	21.8851619590942\\
2.15	21.8851619590942\\
};
\addplot [color=black, line width=1.0pt, only marks, mark=+, mark options={solid, draw=red}, forget plot]
  table[row sep=crcr]{%
1	5.615194452279\\
1	6.77906148336501\\
1	6.78630493450454\\
1	6.78933487401344\\
1	6.94503394207362\\
1	7.00256745549331\\
};
\addplot [color=black, line width=1.0pt, only marks, mark=+, mark options={solid, draw=red}, forget plot]
  table[row sep=crcr]{%
2	19.3362058506939\\
2	19.34281759156\\
2	19.4335641988557\\
2	19.4655334098622\\
2	19.4794454472901\\
2	19.485991065915\\
2	19.6769757139898\\
2	19.6802294554907\\
};
\end{axis}
\end{tikzpicture}%

%% file: BV_sample.tex
%
%
\begin{tikzpicture}

\begin{axis}[%
scale = 0.25,
width=10.0in,
height=3.754in,
at={(1.011in,0.8in)},
scale only axis,
unbounded coords=jump,
xmin=0.5,
xmax=2.5,
xtick={1,2},
xticklabels={{$\gamma_{(a)}$},{$\gamma_{(b)}$}},
xlabel style={font=\color{white!15!black}},
xlabel={SNR=50},
ymin=3.32068333611403,
ymax=56.6755125123286,
axis background/.style={fill=white},
xmajorgrids,
ymajorgrids,
grid style={dotted}
]
\addplot [color=black, dashed, line width=1.0pt, forget plot]
  table[row sep=crcr]{%
1	11.6682943005068\\
1	13.6104588943166\\
};
\addplot [color=black, dashed, line width=1.0pt, forget plot]
  table[row sep=crcr]{%
2	41.6077383113753\\
2	54.2502930043188\\
};
\addplot [color=black, dashed, line width=1.0pt, forget plot]
  table[row sep=crcr]{%
1	7.62594706714906\\
1	9.56439508450616\\
};
\addplot [color=black, dashed, line width=1.0pt, forget plot]
  table[row sep=crcr]{%
2	18.513768837459\\
2	31.499786179299\\
};
\addplot [color=black, line width=1.0pt, forget plot]
  table[row sep=crcr]{%
0.925	13.6104588943166\\
1.075	13.6104588943166\\
};
\addplot [color=black, line width=1.0pt, forget plot]
  table[row sep=crcr]{%
1.925	54.2502930043188\\
2.075	54.2502930043188\\
};
\addplot [color=black, line width=1.0pt, forget plot]
  table[row sep=crcr]{%
0.925	7.62594706714906\\
1.075	7.62594706714906\\
};
\addplot [color=black, line width=1.0pt, forget plot]
  table[row sep=crcr]{%
1.925	18.513768837459\\
2.075	18.513768837459\\
};
\addplot [color=blue, line width=1.0pt, forget plot]
  table[row sep=crcr]{%
0.85	9.56439508450616\\
0.85	11.6682943005068\\
1.15	11.6682943005068\\
1.15	9.56439508450616\\
0.85	9.56439508450616\\
};
\addplot [color=blue, line width=1.0pt, forget plot]
  table[row sep=crcr]{%
1.85	31.499786179299\\
1.85	41.6077383113753\\
2.15	41.6077383113753\\
2.15	31.499786179299\\
1.85	31.499786179299\\
};
\addplot [color=red, line width=1.0pt, forget plot]
  table[row sep=crcr]{%
0.85	11.0245886318406\\
1.15	11.0245886318406\\
};
\addplot [color=red, line width=1.0pt, forget plot]
  table[row sep=crcr]{%
1.85	35.8148691324927\\
2.15	35.8148691324927\\
};
\addplot [color=black, line width=1.0pt, only marks, mark=+, mark options={solid, draw=red}, forget plot]
  table[row sep=crcr]{%
1	5.74590284412379\\
1	15.9784698505641\\
};
\addplot [color=black, line width=1.0pt, only marks, mark=+, mark options={solid, draw=red}, forget plot]
  table[row sep=crcr]{%
nan	nan\\
};
\end{axis}

\end{tikzpicture}%

%% file: BV_approx.tex
%
%
\begin{tikzpicture}

\begin{axis}[%
scale = 0.25,
width=10.0in,
height=3.754in,
at={(1.011in,0.8in)},
scale only axis,
unbounded coords=jump,
xmin=0.5,
xmax=2.5,
xtick={1,2},
xticklabels={{$\gamma_{(a)}$},{$\gamma_{(b)}$}},
xlabel style={font=\color{white!15!black}},
xlabel={SNR=50},
ymin=-0.0796676394136961,
ymax=2.59887897298024,
axis background/.style={fill=white},
xmajorgrids,
ymajorgrids,
grid style={dotted}
]
\addplot [color=black, dashed, line width=1.0pt, forget plot]
  table[row sep=crcr]{%
1	0.0911401995801878\\
1	0.121770840525276\\
};
\addplot [color=black, dashed, line width=1.0pt, forget plot]
  table[row sep=crcr]{%
2	2.02392423095518\\
2	2.47712685423506\\
};
\addplot [color=black, dashed, line width=1.0pt, forget plot]
  table[row sep=crcr]{%
1	0.0420844793314827\\
1	0.0677284314690041\\
};
\addplot [color=black, dashed, line width=1.0pt, forget plot]
  table[row sep=crcr]{%
2	1.43585980394047\\
2	1.72161517191031\\
};
\addplot [color=black, line width=1.0pt, forget plot]
  table[row sep=crcr]{%
0.925	0.121770840525276\\
1.075	0.121770840525276\\
};
\addplot [color=black, line width=1.0pt, forget plot]
  table[row sep=crcr]{%
1.925	2.47712685423506\\
2.075	2.47712685423506\\
};
\addplot [color=black, line width=1.0pt, forget plot]
  table[row sep=crcr]{%
0.925	0.0420844793314827\\
1.075	0.0420844793314827\\
};
\addplot [color=black, line width=1.0pt, forget plot]
  table[row sep=crcr]{%
1.925	1.43585980394047\\
2.075	1.43585980394047\\
};
\addplot [color=blue, line width=1.0pt, forget plot]
  table[row sep=crcr]{%
0.85	0.0677284314690041\\
0.85	0.0911401995801878\\
1.15	0.0911401995801878\\
1.15	0.0677284314690041\\
0.85	0.0677284314690041\\
};
\addplot [color=blue, line width=1.0pt, forget plot]
  table[row sep=crcr]{%
1.85	1.72161517191031\\
1.85	2.02392423095518\\
2.15	2.02392423095518\\
2.15	1.72161517191031\\
1.85	1.72161517191031\\
};
\addplot [color=red, line width=1.0pt, forget plot]
  table[row sep=crcr]{%
0.85	0.0762188891749164\\
1.15	0.0762188891749164\\
};
\addplot [color=red, line width=1.0pt, forget plot]
  table[row sep=crcr]{%
1.85	1.8118318573416\\
2.15	1.8118318573416\\
};
\addplot [color=black, line width=1.0pt, only marks, mark=+, mark options={solid, draw=red}, forget plot]
  table[row sep=crcr]{%
1	0.127604814294436\\
};
\addplot [color=black, line width=1.0pt, only marks, mark=+, mark options={solid, draw=red}, forget plot]
  table[row sep=crcr]{%
nan	nan\\
};
\end{axis}
\end{tikzpicture}%

%% file: analysis_reg.tex
%
%
\begin{tikzpicture}

\begin{axis}[%
scale = 0.25,
width=10.0in,
height=3.754in,
at={(1.011in,0.8in)},
scale only axis,
xmin=0.5,
xmax=3.5,
xtick={1,2,3},
xticklabels={{$\hat{\gamma}$},{$\gamma=0$},{$\gamma_{(b)}$}},
xlabel style={font=\color{white!15!black}},
xlabel={Values of $\gamma$},
ymin=-0.139622294908941,
ymax=2.93221380140608,
axis background/.style={fill=white},
xmajorgrids,
ymajorgrids,
grid style={dotted}
]
\addplot [color=black, dashed, line width=1.0pt, forget plot]
  table[row sep=crcr]{%
1	0.0609458731531023\\
1	0.13516230460768\\
};
\addplot [color=black, dashed, line width=1.0pt, forget plot]
  table[row sep=crcr]{%
2	0.277643867907057\\
2	0.570268890123313\\
};
\addplot [color=black, dashed, line width=1.0pt, forget plot]
  table[row sep=crcr]{%
3	0.0449153307220151\\
3	0.100416674029559\\
};
\addplot [color=black, dashed, line width=1.0pt, forget plot]
  table[row sep=crcr]{%
1	6.61855992399758e-06\\
1	0.00960168943343058\\
};
\addplot [color=black, dashed, line width=1.0pt, forget plot]
  table[row sep=crcr]{%
2	0.00172244770048709\\
2	0.0761156311033156\\
};
\addplot [color=black, dashed, line width=1.0pt, forget plot]
  table[row sep=crcr]{%
3	3.28068186577563e-05\\
3	0.00701137412992677\\
};
\addplot [color=black, line width=1.0pt, forget plot]
  table[row sep=crcr]{%
0.8875	0.13516230460768\\
1.1125	0.13516230460768\\
};
\addplot [color=black, line width=1.0pt, forget plot]
  table[row sep=crcr]{%
1.8875	0.570268890123313\\
2.1125	0.570268890123313\\
};
\addplot [color=black, line width=1.0pt, forget plot]
  table[row sep=crcr]{%
2.8875	0.100416674029559\\
3.1125	0.100416674029559\\
};
\addplot [color=black, line width=1.0pt, forget plot]
  table[row sep=crcr]{%
0.8875	6.61855992399758e-06\\
1.1125	6.61855992399758e-06\\
};
\addplot [color=black, line width=1.0pt, forget plot]
  table[row sep=crcr]{%
1.8875	0.00172244770048709\\
2.1125	0.00172244770048709\\
};
\addplot [color=black, line width=1.0pt, forget plot]
  table[row sep=crcr]{%
2.8875	3.28068186577563e-05\\
3.1125	3.28068186577563e-05\\
};
\addplot [color=blue, line width=1.0pt, forget plot]
  table[row sep=crcr]{%
0.775	0.00960168943343058\\
0.775	0.0609458731531023\\
1.225	0.0609458731531023\\
1.225	0.00960168943343058\\
0.775	0.00960168943343058\\
};
\addplot [color=blue, line width=1.0pt, forget plot]
  table[row sep=crcr]{%
1.775	0.0761156311033156\\
1.775	0.277643867907057\\
2.225	0.277643867907057\\
2.225	0.0761156311033156\\
1.775	0.0761156311033156\\
};
\addplot [color=blue, line width=1.0pt, forget plot]
  table[row sep=crcr]{%
2.775	0.00701137412992677\\
2.775	0.0449153307220151\\
3.225	0.0449153307220151\\
3.225	0.00701137412992677\\
2.775	0.00701137412992677\\
};
\addplot [color=red, line width=1.0pt, forget plot]
  table[row sep=crcr]{%
0.775	0.0281651540571675\\
1.225	0.0281651540571675\\
};
\addplot [color=red, line width=1.0pt, forget plot]
  table[row sep=crcr]{%
1.775	0.158780978019307\\
2.225	0.158780978019307\\
};
\addplot [color=red, line width=1.0pt, forget plot]
  table[row sep=crcr]{%
2.775	0.0193454189564677\\
3.225	0.0193454189564677\\
};
\addplot [color=black, line width=1.0pt, only marks, mark=+, mark options={solid, draw=red}, forget plot]
  table[row sep=crcr]{%
1	0.140617324114121\\
1	0.141028987273406\\
1	0.144146382060493\\
1	0.147330483641218\\
1	0.150582352388758\\
1	0.154423731868446\\
1	0.155289446131162\\
1	0.176629652142073\\
1	0.177333134949568\\
1	0.177724676563134\\
1	0.192013676402483\\
1	0.201440256016574\\
1	0.204336168634067\\
1	0.217988322600188\\
1	0.243997095034279\\
};
\addplot [color=black, line width=1.0pt, only marks, mark=+, mark options={solid, draw=red}, forget plot]
  table[row sep=crcr]{%
2	0.589334688588834\\
2	0.590681269667635\\
2	0.590743172576894\\
2	0.609308750845595\\
2	0.612021274285213\\
2	0.615868443513921\\
2	0.625442553546198\\
2	0.64066428526788\\
2	0.662468920520792\\
2	0.662613875412797\\
2	0.677300203791125\\
2	0.681940041423612\\
2	0.687927591870829\\
2	0.697066377735409\\
2	0.698287872666249\\
2	0.713725960347357\\
2	0.75424390408593\\
2	0.781992270445786\\
2	0.782918635998133\\
2	0.785169725585964\\
2	0.787680895360638\\
2	0.798962992839329\\
2	0.803349275532539\\
2	0.817786915053091\\
2	0.883625644498601\\
2	0.89261568233097\\
2	0.896004866805934\\
2	0.905016435737112\\
2	0.905175204521067\\
2	0.969767704401737\\
2	1.04546708458476\\
2	1.4183340917733\\
2	1.46906416648671\\
2	2.79258488793722\\
};
\addplot [color=black, line width=1.0pt, only marks, mark=+, mark options={solid, draw=red}, forget plot]
  table[row sep=crcr]{%
3	0.102017395375124\\
3	0.102657987510668\\
3	0.102845417873071\\
3	0.103839998518868\\
3	0.10615112482913\\
3	0.1077839562079\\
3	0.112129586287073\\
3	0.123458162566865\\
3	0.123510248054299\\
3	0.123596030283082\\
3	0.124888957888441\\
3	0.128046470035059\\
3	0.13002180490776\\
3	0.138948620217799\\
3	0.143275658470839\\
3	0.14333067103754\\
3	0.14848523160662\\
3	0.195205455172291\\
};
\end{axis}

\end{tikzpicture}%

%% file: main_extended.bbl
\begin{thebibliography}{10}
\providecommand{\url}[1]{#1}
\csname url@samestyle\endcsname
\providecommand{\newblock}{\relax}
\providecommand{\bibinfo}[2]{#2}
\providecommand{\BIBentrySTDinterwordspacing}{\spaceskip=0pt\relax}
\providecommand{\BIBentryALTinterwordstretchfactor}{4}
\providecommand{\BIBentryALTinterwordspacing}{\spaceskip=\fontdimen2\font plus
\BIBentryALTinterwordstretchfactor\fontdimen3\font minus
  \fontdimen4\font\relax}
\providecommand{\BIBforeignlanguage}[2]{{%
\expandafter\ifx\csname l@#1\endcsname\relax
\typeout{** WARNING: IEEEtran.bst: No hyphenation pattern has been}%
\typeout{** loaded for the language `#1'. Using the pattern for}%
\typeout{** the default language instead.}%
\else
\language=\csname l@#1\endcsname
\fi
#2}}
\providecommand{\BIBdecl}{\relax}
\BIBdecl

\bibitem{Pillonetto2014}
G.~Pillonetto, F.~Dinuzzo, T.~Chen, G.~{De Nicolao}, and L.~Ljung, ``Kernel
  methods in system identification, machine learning and function estimation: A
  survey,'' \emph{Automatica}, vol.~50, no.~3, pp. 657--682, 2014.

\bibitem{Hewing2020}
L.~Hewing, K.~P. Wabersich, M.~Menner, and M.~N. Zeilinger, ``Learning-based
  model predictive control: Toward safe learning in control,'' \emph{Annual
  Review of Control, Robotics, and Autonomous Systems}, vol.~3, no.~1, pp.
  269--296, 2020.

\bibitem{TK18}
T.~Koller, F.~Berkenkamp, M.~Turchetta, and A.~Krause, ``Learning-based model
  predictive control for safe exploration,'' in \emph{2018 IEEE Conference on
  Decision and Control (CDC)}, 2018, pp. 6059--6066.

\bibitem{scalableGP}
H.~Liu, Y.-S. Ong, X.~Shen, and J.~Cai, ``When gaussian process meets big data:
  A review of scalable gps,'' \emph{IEEE Transactions on Neural Networks and
  Learning Systems}, vol.~31, no.~11, pp. 4405--4423, 2020.

\bibitem{lazaro2010sparse}
M.~L{\'a}zaro-Gredilla, J.~Qui{\~n}onero-Candela, C.~E. Rasmussen, and A.~R.
  Figueiras-Vidal, ``Sparse spectrum {G}aussian process regression,'' \emph{The
  Journal of Machine Learning Research}, vol.~11, pp. 1865--1881, 2010.

\bibitem{rudi2017}
A.~Rudi and L.~Rosasco, ``Generalization properties of learning with random
  features,'' in \emph{Advances in Neural Information Processing Systems},
  vol.~30.\hskip 1em plus 0.5em minus 0.4em\relax Curran Associates, Inc.,
  2017.

\bibitem{hoerl1970}
A.~E. Hoerl and R.~W. Kennard, ``Ridge regression: Biased estimation for
  nonorthogonal problems,'' \emph{Technometrics}, vol.~12, no.~1, pp. 55--67,
  1970.

\bibitem{Tsybakov2009IntroductionTN}
A.~B. Tsybakov, ``Introduction to nonparametric estimation,'' in \emph{Springer
  series in statistics}, 2009.

\bibitem{pmlr-v70-pan17a}
Y.~Pan, X.~Yan, E.~A. Theodorou, and B.~Boots, ``Prediction under uncertainty
  in sparse spectrum {G}aussian processes with applications to filtering and
  control,'' in \emph{Proceedings of the 34th International Conference on
  Machine Learning}, ser. PMLR, vol.~70, August 2017, pp. 2760--2768.

\bibitem{arcari2021}
E.~Arcari, A.~Scampicchio, A.~Carron, and M.~N. Zeilinger, ``Bayesian
  multi-task learning using finite-dimensional models: A comparative study,''
  in \emph{2021 60th IEEE Conference on Decision and Control (CDC)}, 2021, pp.
  2218--2225.

\bibitem{arcari22robot}
\BIBentryALTinterwordspacing
E.~Arcari, M.~V. Minniti, A.~Scampicchio, A.~Carron, F.~Farshidian, M.~Hutter,
  and M.~N. Zeilinger, ``Bayesian multi-task learning mpc for robotic mobile
  manipulation,'' 2022. [Online]. Available:
  \url{https://arxiv.org/abs/2211.10270}
\BIBentrySTDinterwordspacing

\bibitem{Cucker02}
F.~Cucker and S.~Smale, ``On the mathematical foundations of learning,''
  \emph{Bulletin of the American Mathematical Society}, vol.~39, pp. 1--49,
  2002.

\bibitem{cuckerzhou2007}
F.~Cucker and D.~X. Zhou, \emph{Learning Theory: An Approximation Theory
  Viewpoint}, ser. Cambridge Monographs on Applied and Computational
  Mathematics.\hskip 1em plus 0.5em minus 0.4em\relax Cambridge University
  Press, 2007.

\bibitem{Smale2007LearningTE}
S.~Smale and D.-X. Zhou, ``Learning theory estimates via integral operators and
  their approximations,'' \emph{Constructive Approximation}, vol.~26, pp.
  153--172, 2007.

\bibitem{lunardi2009}
A.~Lunardi, \emph{Interpolation Theory}, ser. Publications of the {S}cuola
  {N}ormale di {P}isa.\hskip 1em plus 0.5em minus 0.4em\relax Edizioni della
  Normale Pisa, 2009.

\bibitem{Boucheron2004}
S.~Boucheron, G.~Lugosi, and O.~Bousquet, \emph{Concentration
  Inequalities}.\hskip 1em plus 0.5em minus 0.4em\relax Berlin, Heidelberg:
  Springer Berlin Heidelberg, 2004, pp. 208--240.

\bibitem{Wu06learningrates}
Q.~Wu, Y.~Ying, and D.-X. Zhou, ``Learning rates of least-square regularized
  regression,'' \emph{Foundations of Computational Mathematics}, pp. 171--192,
  2006.

\bibitem{cucker2008bestchoices}
F.~Cucker and S.~Smale, ``Best choices for regularization parameters in
  learning theory: On the bias—variance problem,'' \emph{Foundations of
  Computational Mathematics}, vol.~2, pp. 413--428, March 2008.

\bibitem{mendelson2010}
S.~Mendelson and J.~Neeman, ``{Regularization in kernel learning},'' \emph{The
  Annals of Statistics}, vol.~38, no.~1, pp. 526 -- 565, 2010.

\bibitem{WANG201155}
C.~Wang and D.-X. Zhou, ``Optimal learning rates for least squares regularized
  regression with unbounded sampling,'' \emph{Journal of Complexity}, vol.~27,
  no.~1, pp. 55--67, 2011.

\bibitem{Caponnetto2007OptimalRF}
A.~Caponnetto and E.~de~Vito, ``Optimal rates for the regularized least-squares
  algorithm,'' \emph{Foundations of Computational Mathematics}, vol.~7, pp.
  331--368, 2007.

\bibitem{Wang2011}
C.~Wang and D.-X. Zhou, ``\BIBforeignlanguage{English}{Optimal learning rates
  for least squares regularized regression with unbounded sampling},''
  \emph{\BIBforeignlanguage{English}{Journal of Complexity}}, vol.~27, no.~1,
  pp. 55--67, 2011.

\bibitem{Guo2013ConcentrationEF}
Z.-C. Guo and D.-X. Zhou, ``Concentration estimates for learning with unbounded
  sampling,'' \emph{Advances in Computational Mathematics}, vol.~38, pp.
  207--223, 2013.

\bibitem{zhou2003}
D.-X. Zhou, ``Capacity of reproducing kernel spaces in learning theory,''
  \emph{IEEE Transactions on Information Theory}, vol.~49, pp. 1743 -- 1752,
  August 2003.

\bibitem{akhiezer2013theory}
N.~Akhiezer and I.~Glazman, \emph{Theory of Linear Operators in Hilbert Space},
  ser. Dover Books on Mathematics.\hskip 1em plus 0.5em minus 0.4em\relax Dover
  Publications, 2013.

\bibitem{aronszajn50reproducing}
N.~Aronszajn, ``Theory of reproducing kernels,'' \emph{Transactions of the
  American Mathematical Society}, vol.~68, no.~3, pp. 337--404, 1950.

\bibitem{steinwartSVM}
I.~Steinwart and A.~Christmann, \emph{Support Vector Machines}, 1st~ed.\hskip
  1em plus 0.5em minus 0.4em\relax Springer Publishing Company, Incorporated,
  2008.

\bibitem{Lin2017distributed}
S.-B. Lin, X.~Guo, and D.-X. Zhou, ``Distributed learning with regularized
  least squares,'' \emph{Journal of Machine Learning Research}, vol.~18,
  no.~92, pp. 1--31, 2017.

\bibitem{Davies2007NonparametricRC}
P.~L. Davies, A.~Kovac, and M.~Meise, ``Nonparametric regression, confidence
  regions and regularization,'' \emph{Annals of Statistics}, vol.~37, pp.
  2597--2625, 2007.

\bibitem{Niyogi99}
P.~Niyogi and F.~Girosi, ``Generalization bounds for function approximation
  from scattered noisy data,'' \emph{Advances in Computational Mathematics},
  vol.~10, no.~1, 1999.

\bibitem{berberian1961introduction}
S.~Berberian, \emph{Introduction to Hilbert Space}.\hskip 1em plus 0.5em minus
  0.4em\relax Oxford University Press, 1961.

\bibitem{wahba2019representer}
G.~Wahba and Y.~Wang, \emph{Representer Theorem}.\hskip 1em plus 0.5em minus
  0.4em\relax American Cancer Society, 2019, pp. 1--11.

\bibitem{smaleshannonII}
S.~Smale and D.-X. Zhou, ``Shannon sampling {II}: Connections to learning
  theory,'' \emph{Applied and Computational Harmonic Analysis}, vol.~19, no.~3,
  pp. 285--302, 2005.

\bibitem{smaleshannonI}
------, ``Shannon sampling and function reconstruction from point values,''
  \emph{Bulletin of The American Mathematical Society}, vol.~41, pp. 279--306,
  July 2004.

\bibitem{zhang2005effdim}
T.~Zhang, ``{Learning Bounds for Kernel Regression Using Effective Data
  Dimensionality},'' \emph{Neural Computation}, vol.~17, no.~9, pp. 2077--2098,
  September 2005.

\bibitem{MU2018381}
B.~Mu, T.~Chen, and L.~Ljung, ``On asymptotic properties of hyperparameter
  estimators for kernel-based regularization methods,'' \emph{Automatica},
  vol.~94, pp. 381--395, 2018.

\bibitem{gilks1995markov}
W.~Gilks, S.~Richardson, and D.~Spiegelhalter, \emph{Markov Chain Monte Carlo
  in Practice}, ser. Chapman \& Hall/CRC Interdisciplinary Statistics.\hskip
  1em plus 0.5em minus 0.4em\relax Taylor \& Francis, 1995.

\bibitem{PILLONETTO2015106}
G.~Pillonetto and A.~Chiuso, ``Tuning complexity in regularized kernel-based
  regression and linear system identification: The robustness of the marginal
  likelihood estimator,'' \emph{Automatica}, vol.~58, pp. 106--117, 2015.

\bibitem{maritz2018empirical}
J.~Maritz, \emph{Empirical Bayes Methods with Applications}.\hskip 1em plus
  0.5em minus 0.4em\relax CRC Press, 2018.

\end{thebibliography}
